\documentclass[11pt,reqno]{amsproc}
\allowdisplaybreaks
\usepackage{color}
\usepackage{colortbl}
\usepackage{fullpage}
\usepackage{graphicx}
\usepackage{textcomp}
\usepackage{subfigure}
\usepackage{enumerate}
\usepackage{algorithm2e}
\usepackage[debug=false, colorlinks=true, pdfstartview=FitV,
linkcolor=blue, citecolor=blue, urlcolor=blue]{hyperref}
\usepackage[semicolon,square,authoryear,sort]{natbib}

\usepackage[usenames,dvipsnames,table]{xcolor}
\usepackage[most]{tcolorbox}

\usepackage{multirow} %for multirow
\usepackage{rotating} %for sideways, rotating "samples" at 90°"
\usepackage{bigstrut} %for bigstrut
\usepackage{hhline} 

\usepackage[doipre={DOI:~}]{uri}

\usepackage{morefloats}

\newlength{\drop}
\definecolor{amethyst}{rgb}{0.6, 0.4, 0.8}
\definecolor{burgundy}{rgb}{0.5, 0.0, 0.13}

\title{\textbf{How to pose material design problems for flow through porous media applications?: Sensitivity of dissipation rate to medium's permeability holds the key}}

\author{\textbf{Kalyana B.~Nakshatrala} \\
  {\small Associate Professor, Department of Civil and Environmental Engineering \\
  University of Houston, Houston, Texas 77204, USA.}\\
  {\small email:~\texttt{knakshatrala@uh.edu}, phone: +1-713-743-4418}}


\begin{document}

%===========================;
%  Title page of the paper  ;
%===========================;
\begin{titlepage}
  \drop=0.1\textheight
  \centering
  \vspace*{\baselineskip}
  \rule{\textwidth}{1.6pt}\vspace*{-\baselineskip}\vspace*{2pt}
  \rule{\textwidth}{0.4pt}\\[\baselineskip]
      {\Large \textbf{\color{burgundy}
      How to pose material design problems for flow through \\[0.3\baselineskip] porous media applications?: Sensitivity of dissipation rate \\[0.3\baselineskip] 
      to medium's permeability holds the key}}\\[0.3\baselineskip]
      \rule{\textwidth}{0.4pt}\vspace*{-\baselineskip}\vspace{3.2pt}
      \rule{\textwidth}{1.6pt}\\[\baselineskip]
      \scshape
      An e-print of this paper is available on arXiv:~2110.10827. \par
      \vspace*{1\baselineskip}
      Authored by \\[\baselineskip]

  {\Large K.~B.~Nakshatrala\par}
  {\itshape Department of Civil \& Environmental Engineering \\
  University of Houston, Houston, Texas 77204. \\
  \textbf{phone:} +1-713-743-4418, \textbf{e-mail:} knakshatrala@uh.edu \\
  \textbf{website:} http://www.cive.uh.edu/faculty/nakshatrala\\
  \textbf{ORCID ID}: https://orcid.org/0000-0002-1744-6338}\\[\baselineskip]
  
  \vfill
  {\scshape 2021} \\
  {\small Computational \& Applied Mechanics Laboratory} \par
\end{titlepage}

%=========================;
%  Abstract of the paper  ;
%=========================;
\begin{abstract}
Recent studies have advocated using the total dissipation rate under topology optimization to realize material designs involving the flow of fluids through porous media. However, these studies decided how to pose the design problem, such as maximizing the total dissipation rate for some situations while minimizing for others, by solving one-dimensional problems and justifying their choices using numerical experiments. The rigor is lacking---a bottleneck for further scientific advancements to computational material design. This paper provides the missing theoretical justification. We identify four classes of boundary value problems using the adjoint state method and analytically calculate the sensitivity of the total dissipation rate to the permeability field. For two of those classes in which the flow of fluids is pressure-driven, the sensitivity is positive---the total dissipation rate increases if the medium's permeability increases. While for the other two classes, in which the flow is velocity-driven, the trend is the opposite. These sensitivities provide rigorous answers to the central question: how to pose a material design problem for flow through porous media applications. The impact of our work is multi-fold. \emph{First}, this study further elevates the role of the dissipation rate in posing well-posed material design problems using topology optimization. \emph{Second}, besides the theoretical significance, the results benefit computational scientists and practitioners to realize optimal designs. \emph{Third}, given their simplicity yet far-reaching impact, both the approach and results possess immense pedagogical value.
\end{abstract}

\maketitle

\vspace{-0.45in}

%==================================;
%  Include all the sections below  ;
%==================================;

%**********************************************;
%                                              ;
%  NAME                                        ;
%    S1_Posing_Intro.tex                       ;
%                                              ;
%  WRITTEN BY                                  ;
%    Kalyana B. Nakshatrala                    ;
%                                              ;
%**********************************************;
\section{INTRODUCTION AND MOTIVATION}
\label{Sec:S1_Posing_Intro}

The flow of fluids through porous media is ubiquitous in both natural and synthetic worlds. Examples include nutrient transport in tissues \citep{friedman2008principles,dey2021mathematical}, the flow of hydrocarbons in a petroleum reservoir \citep{ewing1983mathematics}, water filtration using membranes \citep{chen2021flow,herterich2015tailoring}, and workings of facial masks to protect against transmission of COVID-19 and air pollution \citep{verma2020visualizing,krishan2021efficacy}, to name a few. Given the field’s wide-ranging importance, researchers have been developing mathematical models \citep{coussy2004poromechanics,Rajagopal_2007,de2012theory}, numerical formulations \citep{brezzi1991mixed,masud2002stabilized,chang2017modification}, computer packages \citep{COMSOL}, and \emph{a posteriori} accuracy measures \citep{shabouei_nakshatrala_cicp} to study such problems. Material design is always a primary focus given the constant quest to realize new functionalities and optimize associated devices’ performance.

However, ``optimal" design will have different connotations depending on the application and the functionality one wants to realize. For instance, in water filtration using membranes, the surface of the pores is coated with (chemical) functional groups so that microbes, which have surface charges, are removed via electrostatic attraction or repulsion. The question in such applications is how to manipulate pore surfaces to extract pathogens optimally. In applications involving thermal regulation (e.g., heat exchangers), researchers alter the multi-phase characteristics of the fluids to maximize cooling efficiency. Yet, in some other applications, one designs the very geometrical nature of the pore structure. In this paper, the material design optimizes the dissipation rate via the spatial distribution of porous materials with different permeabilities.

In the past, the limitations of fabrication techniques did not allow researchers to take full advantage of novel designs, which often have complex microstructures. With recent advances in additive manufacturing, the fabrication of complex layouts is now possible \citep{sigmund2009manufacturing,liu2016survey}. So, now that manufacturing is not the bottleneck, researchers look for intelligent designs: realized via designers' ingenuity or mathematical-driven solution procedures.

Topology optimization is a leading mathematical framework for designing material systems \citep{bendsoe_sigmund_2013topology}. Although the early studies focused on structural optimization, topology optimization is increasingly popular in several other fields involving the flow of fluids through porous media; these emerging fields include the design of microfluidic devices \citep{kreissl2011topology,alexandersen2020review} and medical devices 
such as dialyzer (filter) used for hemodialysis \citep{alonso2021topology,cancilla2022porous}. However, many works apply topology optimization in these relatively new application areas by mimicking the material design problem from solid mechanics. For example, some of these works replace the compliance with a similar-looking quadratic functional or replace the \emph{principle of minimum potential energy} from structural mechanics with the \emph{principle of minimum power} \citep{guest2006topology,phatak2021optimal}. However, such a translation from solid mechanics to flow through porous media may not be physically meaningful and sometimes not even consistent with the state equations. For example, the minimum power principle does not apply to nonlinear porous models but is limited to linear models such as Darcy and Darcy-Brinkman. More importantly, the said principle is not a physical law but rather a restatement of the governing equations for linear models using the calculus of variations \citep{phatak2021optimal}. 

In pursuit of an objective function to drive the material design for flow through porous media applications, \cite{phatak2021optimal} have promoted the \emph{total dissipation rate} as an ideal candidate to pose material design problems using topology optimization. This fundamental physical quantity has a universal appeal---applicable even to nonlinear models---which is not the case with the alternatives such as the principle of minimum power. Because of its recent popularity, the total dissipation rate still poses several unanswered questions for defining well-founded material design problems. Three such questions are:
\begin{enumerate}
\item[(a)] whether should we maximize or minimize the total dissipation rate, 
\item[(b)] should the volumetric bound constraint be on the use of high- or low-permeability material, and
\item[(c)] when does the optimization problem render nontrivial designs?
\end{enumerate}

The cited paper addressed the first two design questions by considering one-dimensional problems under the Darcy model and supported by numerical results on more complicated boundary value problems. Alluding to the success under the Darcy model, \citet{phatak2021effect} have adopted a similar way of posing the material design for the Darcy-Brinkman model. Although the one-dimensional problems offered sound guidance and numerical results provided the proof of concept, a lack of rigor hinders the progress on theoretical and practical fronts.  

Thus, the central aim of this paper is to fill this lacuna by providing an in-depth analysis of \textit{how to pose the material design problem for flow through porous media applications}. We consider Darcy and Darcy-Brinkman models---the two prominent mathematical models describing the flow of incompressible fluids through porous media. Our approach first estimates the \emph{design sensitivity}---the sensitivity of the total dissipation rate (i.e., the objective function in a material design problem) to the permeability field (i.e., the design variable)---using the \emph{adjoint state method}\footnote{The adjoint state method is also referred to as the adjoint variable method or simply the adjoint method.}.

The adjoint state method has an eventful history, with its first application traced to the \emph{optimal control of systems governed by partial differential equations}---a seminal work by \citet{lions1971optimal}. Since then, researchers have applied this method in the design of structural systems \citep{haug1986design}, seismology \citep{tromp2005seismic}, topology optimization \citep{nakshatrala2013nonlinear,nakshatrala2016nonlinear,jensen2014consistency}, reliability studies and parameter identification \citep{tortorelli1994design}, check-pointing in computer codes \citep{wang2009minimal}, machine learning \citep{pmlr-v119-zhuang20a}, PDE-constrained optimization \citep{hinze2008optimization,bradley2013pde}, and inverse problems in subsurface modeling \citep{sun2013inverse}. The main advantage of the adjoint state method is in the calculation of the sensitivity of an objective function to a set of design parameters. The method circumvents the explicit calculation of the sensitivity of the solution fields by introducing Lagrange multipliers, calculated by solving an  adjoint problem (also referred to as the adjoint state problem), which itself is a system of partial differential equations. 

Our analysis---based on the adjoint state method---identifies four principal classes of boundary problems. Two of them are pressure-driven, and the other two are velocity-driven; the adjoint problem is analytically solvable for these classes of problems. As shown in this paper, the answers to the questions raised above---how the total dissipation rate alters with the permeability field and how to pose material design problems---depend on the class of problems. Notably, for the two pressure-driven classes of problems, the total dissipation rate increases with an increase in the permeability field; in comparison, the reverse trend holds for the two velocity-driven classes of problems. Thus, for example, trivial designs occur under different scenarios for pressure-driven and velocity-driven problems. These findings have theoretical, numerical, practical, and pedagogical significance. The research will place the total dissipation rate for material design problems on a sound theoretical footing, achieved by rigorously answering the aforementioned design-related questions.

The plan for the rest of this paper is as follows. Section \ref{Sec:S2_Posing_Darcy_equations} presents Darcy equations and the dissipation rate under the Darcy model. Section \ref{Sec:S3_Posing_Four_classes} identifies four classes of boundary value problems, which are central to much of the discussion in this paper. Using the adjoint state method, Section \ref{Sec:S4_Posing_Adjoint_method} derives the sensitivity of the total dissipation to the permeability field and solves the adjoint problem analytically for the four classes of boundary problems. Section \ref{Sec:S5_Posing_Brinkman} shows that these findings translate even to the Darcy-Brinkman model. By synthesizing the results from the previous sections, Section \ref{Sec:S6_Posing_Material_design} addresses the central and related design questions with mathematical justifications. Section \ref{Sec:S7_Posing_CR} offers concluding remarks.

%**********************************************;
%                                              ;
%  NAME                                        ;
%    S2_Posing_Darcy_equations.tex             ;
%                                              ;
%  WRITTEN BY                                  ;
%    Kalyana B. Nakshatrala                    ;
%                                              ;
%**********************************************;
\section{DARCY EQUATIONS AND TOTAL DISSIPATION RATE}
\label{Sec:S2_Posing_Darcy_equations}
Consider the flow of an incompressible single-phase fluid in a porous domain $\Omega$ with boundary $\partial \Omega$. The fluid's density and coefficient of viscosity are denoted as $\rho$ and $\mu$, respectively. $k(\mathbf{x})$ denotes the permeability of the porous domain at a spatial point $\mathbf{x} \in \Omega$. The spatial divergence and gradient operators are denoted by $\mathrm{div}[\cdot]$ and $\mathrm{grad}[\cdot]$, respectively. 

A general forward problem\footnote{Alternative names used in the literature for the forward problem are the direct problem, state equations, and analysis problem.} under the Darcy model---the so-called Darcy equations---reads: Find velocity vector field $\mathbf{v}$ and pressure scalar field $p$ satisfying 
%-----------------------------;
%  Equation: Darcy equations  ;
%-----------------------------;
\begin{subequations}
\begin{alignat}{2}
    \label{Eqn:Posing_Darcy_BoLM} 
    &\frac{\mu}{k(\mathbf{x})} \mathbf{v} + \mathrm{grad}[p] = \rho \mathbf{b}(\mathbf{x}) 
    &&\quad \mathrm{in} \; \Omega \\
    \label{Eqn:Posing_Darcy_Continuity} 
    &\mathrm{div}[\mathbf{v}] = 0 
    &&\quad \mathrm{in} \; \Omega \\
    \label{Eqn:Posing_Darcy_pBC} 
    &p = p^{\mathrm{p}}(\mathbf{x}) 
    &&\quad \mathrm{on} \; \Gamma^{p} \\
    \label{Eqn:Posing_Darcy_vBC} 
    &\mathbf{v}\cdot \widehat{\mathbf{n}}(\mathbf{x}) = v_n^{\mathrm{p}}(\mathbf{x}) 
    &&\quad \mathrm{on} \; \Gamma^{v} 
\end{alignat}
\end{subequations}
where $\mathbf{b}(\mathbf{x})$ is the body force vector field; $p^{\mathrm{p}}(\mathbf{x})$ and $v_n^{\mathrm{p}}(\mathbf{x})$ are, respectively, the prescribed pressure and normal component of the velocity on the boundary; $\widehat{\mathbf{n}}(\mathbf{x})$ is the outward unit vector on the boundary; 
and $\Gamma^{p}$ and $\Gamma^{v}$ are the parts of boundary on which pressure and velocity boundary conditions are prescribed, respectively. For mathematical well posedness, the two parts of the boundary satisfy: 
%-------------------------------------------;
%  Equation: Decomposition of the boundary  ;
%-------------------------------------------;
\begin{align}
    \Gamma^{p} \cup \Gamma^{v} = \partial \Omega 
    \quad \mathrm{and} \quad 
    \Gamma^{p} \cap \Gamma^{v} = \emptyset 
\end{align}

Solutions to Darcy equations are \emph{unique}, except for the case when the velocity is prescribed on the entire boundary (i.e., $\Gamma^{v} = \partial \Omega$). For such a case, a solution exists \emph{only} when the following compatibility condition is met:
%-------------------------------------;
%  Equation: Compatibility condition  ;
%-------------------------------------;
\begin{align}
    \int_{\Gamma^{v} = \partial \Omega} v_n^{\mathrm{p}}(\mathbf{x}) \, \mathrm{d}\Gamma = 0 
\end{align}
The rationale behind the above compatibility condition is the divergence theorem. To wit, 
%------------------------------------------------;
%  Equation: Justifying compatibility condition  ;
%------------------------------------------------;
\begin{align}
    0 = \int_{\Omega} \mathrm{div}[\mathbf{v}]\,\mathrm{d} \Omega = \int_{\partial \Omega} \mathbf{v} \cdot \widehat{\mathbf{n}}(\mathbf{x}) \, \mathrm{d} \Gamma = \int_{\Gamma^{v} = \partial \Omega} v_n^{\mathrm{p}}(\mathbf{x}) \, 
    \mathrm{d} \Gamma 
\end{align}
For the case $\Gamma^{v} = \partial \Omega$ with a  compatible velocity boundary condition, the pressure field can be uniquely found only up to a constant; however, the velocity solution field is still unique. The uniqueness of solutions will be handy in ascertaining the solution for adjoint variables. 

In a forward problem, one solves Darcy equations to get the solution (i.e., velocity and pressure) fields for a fixed set of inputs, including a prescribed permeability field. For such analyses (with assigned permeability fields), solution fields are mere spatial fields. However, in a material design problem, say topology optimization, the permeability field is the design variable---an unknown by itself. During design iterations, the permeability field changes, altering the solution fields in each iteration. Regarding this paper, we address changes to the total dissipation rate upon changes to the permeability field. To quantify such changes (i.e., to calculate sensitivities), one must account for the implicit dependence of the dissipation rate on the permeability field via the solution fields. Said differently, one cannot treat the solution fields as mere spatial fields, depending only on the spatial location; a change to the permeability field alters the solution fields.

Hence, we distinguish the solution fields from arbitrary spatial fields to emphasize their (implicit) dependence on the permeability field. From hereon, we denote the solutions fields---the Darcy velocity and the associated pressure---as $\mathbf{v}_{\mathrm{D}}\big(\mathbf{x};k(\mathbf{x})\big)$ and $p_{\mathrm{D}}\big(\mathbf{x};k(\mathbf{x})\big)$. Note the semicolon notation: the solution fields depend explicitly on the quantities on the left side of the semicolon, while the quantity on the right side of the semicolon is the design variable---the permeability field in this paper.

For given permeability $k(\mathbf{x})$ and (arbitrary) vector $\mathbf{v}(\mathbf{x})$ fields, the rate of dissipation under the  Darcy model takes the following mathematical form:
%--------------------------------------;
%  Equation: Dissipation rate density  ;
%--------------------------------------;
\begin{align}
\varphi_{\mathrm{D}} \big(\mathbf{x},k(\mathbf{x}),\mathbf{v}(\mathbf{x})\big) 
= \frac{\mu}{k(\mathbf{x})} \mathbf{v}(\mathbf{x}) \cdot \mathbf{v}(\mathbf{x})
\end{align}
Given a dissipation rate density, the corresponding total dissipation rate, calculated over the entire domain, reads: 
%-----------------------------------------------;
%  Equation: Total dissipation rate functional  ;
%-----------------------------------------------;
\begin{align}
 \Phi_{\mathrm{D}} \big(k(\mathbf{x}),\mathbf{v}(\mathbf{x})\big)
:= \int_{\Omega} \varphi_{\mathrm{D}} \big(\mathbf{x},k(\mathbf{x}),\mathbf{v}(\mathbf{x})\big) \, \mathrm{d} \Omega 
= \int_{\Omega} \frac{\mu}{k(\mathbf{x})} \mathbf{v}(\mathbf{x}) \cdot \mathbf{v}(\mathbf{x}) \, \mathrm{d} \Omega 
\end{align}

We introduce $\Phi_{\mathrm{D}}[k(\mathbf{x})]$---referred to as the total dissipation rate \emph{functional}\footnote{A functional is a function of functions; see \citep{gelfand2000calculus}.}---to denote the total dissipation rate expended under the \emph{Darcy model} by the \emph{Darcy velocity} field with a permeability field $k(\mathbf{x})$. Mathematically,  
%-----------------------------------------------;
%  Equation: Total dissipation rate functional  ;
%-----------------------------------------------;
\begin{align}
\Phi_{\mathrm{D}}\big[k(\mathbf{x})\big]
:= \Phi_{\mathrm{D}} \Big(k(\mathbf{x}),\mathbf{v}_{\mathrm{D}}\big(\mathbf{x};k(\mathbf{x})\big)\Big) 
= \int_{\Omega} \frac{\mu}{k(\mathbf{x})} \mathbf{v}_{\mathrm{D}}\big(\mathbf{x};k(\mathbf{x})\big) \cdot \mathbf{v}_{\mathrm{D}}\big(\mathbf{x};k(\mathbf{x})\big) \, \mathrm{d} \Omega 
\end{align}

%======================================================;
%  Subsection: A notation for calculus of functionals  ;
%======================================================;
\subsection{A notation for calculus of functionals} 
The (Fr\'echet)  derivative of a functional  $\Phi[k(\mathbf{x})]$ at $k(\mathbf{x})$ is defined as \citep{spivak2018calculus}: 
\begin{align}
    \lim_{\|\Delta k(\mathbf{x})\| \rightarrow 0} \frac{\Phi[k(\mathbf{x}) + \Delta k(\mathbf{x})] 
    - \Phi[k(\mathbf{x})] 
    - D\Phi[k(\mathbf{x})] \cdot \Delta k (\mathbf{x})}{\|\Delta k(\mathbf{x})\|} = 0
\end{align}
for any perturbation field $\Delta k(\mathbf{x})$, and $\|\cdot\|$ is the standard Euclidean norm. For a given scalar field of the form $p(\mathbf{x};k(\mathbf{x}))$ (e.g., the pressure field satisfying the forward problem), the derivatives $D_1p$ and $D_2p$ are respectively defined as follows:
\begin{align}
    \label{Eqn:Posing_D1_definition}
    &\lim_{\|\Delta \mathbf{x}\| \rightarrow 0}  \frac{p\big(\mathbf{x} + \Delta \mathbf{x};k(\mathbf{x})\big) 
   - p\big(\mathbf{x};k(\mathbf{x})\big)  
    - D_{1}p\big(\mathbf{x};k(\mathbf{x})\big)  
    \cdot \Delta \mathbf{x}}{\|\Delta \mathbf{x}\|} = 0 \\
    &\lim_{\|\Delta k(\mathbf{x})\| \rightarrow 0}  \frac{p\big(\mathbf{x};k(\mathbf{x}) + \Delta k(\mathbf{x})\big) 
   - p\big(\mathbf{x};k(\mathbf{x})\big)  
    - D_{2}p\big(\mathbf{x};k(\mathbf{x})\big)  
    \cdot \Delta k(\mathbf{x})}{\|\Delta k(\mathbf{x})\|} = 0
\end{align}
For convenience, we denote these two derivatives as follows: 
\begin{align}
&\mathrm{grad}[p(\mathbf{x};k(\mathbf{x}))] \equiv D_{1}p(\mathbf{x};k(\mathbf{x})) \\
&p^{\prime}(\mathbf{x};k(\mathbf{x})) \equiv D_{2}p(\mathbf{x};k(\mathbf{x})) 
\end{align}    

Likewise, one can define $\mathrm{grad}[\mathbf{v}(\mathbf{x};k(\mathbf{x}))] = D_1\mathbf{v}(\mathbf{x};k(\mathbf{x}))$ and $\mathbf{v}^{\prime}(\mathbf{x};k(\mathbf{x})) = D_2 \mathbf{v}(\mathbf{x};k(\mathbf{x}))$ for a vector field $\mathbf{v}(\mathbf{x};k(\mathbf{x}))$. As done in tensor algebra, the divergence can then be defined as follows: 
\begin{align}
    \mathrm{div}[\mathbf{v}] = \mathrm{tr}\big[\mathrm{grad}[\mathbf{v}]\big]
\end{align}
where $\mathrm{tr}[\cdot]$ is the trace operator for second-order tensors \citep{chadwick2012continuum}.

For a given scalar field of the form  $\varphi\Big(\mathbf{x},k(\mathbf{x}),\mathbf{v}(\mathbf{x};k(\mathbf{x})\big)\Big)$, we define the derivatives $D_2\varphi$ and $D_3\varphi$ as follows:
{\small
\begin{align}
    &\lim_{\|\Delta k(\mathbf{x})\| \rightarrow 0}  \frac{\varphi\big(\mathbf{x},k(\mathbf{x})+\Delta k,\mathbf{v}\big(\mathbf{x};k(\mathbf{x})\big)\Big)
   - \varphi\Big(\mathbf{x},k(\mathbf{x}),\mathbf{v}\big(\mathbf{x};k(\mathbf{x})\big)\Big) 
    - D_{2}\varphi\Big(\mathbf{x},k(\mathbf{x}),\mathbf{v}\big(\mathbf{x};k(\mathbf{x})\big)\Big)  
    \cdot \Delta k(\mathbf{x})}{\|\Delta k(\mathbf{x})\|} = 0 \\
    &\lim_{\|\Delta \mathbf{v}(\mathbf{x})\| \rightarrow 0}  \frac{\varphi\Big(\mathbf{x},k(\mathbf{x}),\mathbf{v}(\mathbf{x};k(\mathbf{x})) + \Delta \mathbf{v}(\mathbf{x})\Big) 
   - \varphi\Big(\mathbf{x},k(\mathbf{x}),\mathbf{v}\big(\mathbf{x};k(\mathbf{x})\big)\Big)
    - D_{2}\varphi\Big(\mathbf{x},k(\mathbf{x}),\mathbf{v}\big(\mathbf{x};k(\mathbf{x})\big) \Big)
    \cdot \Delta \mathbf{v}(\mathbf{x})}{\|\Delta \mathbf{v}(\mathbf{x})\|} = 0
\end{align}
}
The definition for $D_1\varphi$ is similar to the one provided before (cf.~Eq.~\eqref{Eqn:Posing_D1_definition}), and $\mathrm{grad}[\varphi]\equiv D_1\varphi$. Appealing to the chain rule, we define the total derivative of $\varphi$ with respect to the permeability field as follows: 
%------------------------------;
%  Equation: Total dependence  ;
%------------------------------;
{\small
\begin{align}
&\varphi^{\prime}\Big(\mathbf{x},k(\mathbf{x}),\mathbf{v}\big(\mathbf{x};k(\mathbf{x})\big)\Big)
:= D_{2}\varphi\Big(\mathbf{x},k(\mathbf{x}),\mathbf{v}\big(\mathbf{x};k(\mathbf{x})\big)\Big)
+ D_{3}\varphi\Big(\mathbf{x},k(\mathbf{x}),\mathbf{v}\big(\mathbf{x};k(\mathbf{x})\big)\Big) \cdot 
D_{2}\mathbf{v}\big(\mathbf{x};k(\mathbf{x})\big)
\end{align}}
The first term on the right side of the above equation accounts for explicit dependence, while the second term accounts for implicit dependence.

%**********************************************;
%                                              ;
%  NAME                                        ;
%    S3_Posing_Four_classes.tex                ;
%                                              ;
%  WRITTEN BY                                  ;
%    Kalyana B. Nakshatrala                    ;
%                                              ;
%**********************************************;
\section{FOUR CLASSES OF BOUNDARY VALUE PROBLEMS}
\label{Sec:S3_Posing_Four_classes}
For further analysis, we identify the following four classes of boundary value problems under the Darcy model:\footnote{The rationale behind this classification will be apparent in subsequent sections, where we solve the adjoint problem (e.g., see \S\ref{Sec:S4_Posing_Adjoint_method}).}
%-------------------------------------;
%  Enumerate: Classification of BVPs  ;
%-------------------------------------;
\begin{enumerate}
\item \textbf{Class A:} Prescription of a pressure loading on the entire boundary. Stated mathematically, 
\begin{align}
\Gamma^{v} = \emptyset
\end{align}
\item \textbf{Class B:} Prescription of pressure loading on a part of the boundary while homogeneous velocity boundary conditions on the rest. That is, 
\begin{align}
\Gamma^{p} \neq \emptyset, \; 
\Gamma^{v} \neq \emptyset, 
\; \mathrm{and} \;  v_n^{\mathrm{p}}(\mathbf{x}) = 0 \; \mathrm{on} \;  \Gamma^{v}
\end{align}
\item \textbf{Class C:} Prescription of a compatible velocity boundary condition on the entire boundary and requiring the prescribed body force to be a conservative vector field. Mathematically, 
\begin{align}
\Gamma^{p} = \emptyset, \; \int_{\Gamma^{v} = \partial \Omega} v_n^{\mathrm{p}}(\mathbf{x}) \, \mathrm{d} \Gamma = 0, \; \mathrm{and} \; \rho\mathbf{b}(\mathbf{x}) = - \mathrm{grad}[\psi] 
\end{align}
where $\psi(\mathbf{x})$ is a scalar field. 
\item \textbf{Class D:} Prescription of velocity boundary conditions on a part of the boundary with zero pressure loading on the rest, and requiring the body force is zero. Written mathematically, 
\begin{align}
\Gamma^{p} \neq \emptyset, \; \Gamma^{v} \neq 0, \;  p^{\mathrm{p}}(\mathbf{x}) = 0 \; \mathrm{on} \; \Gamma^{p}, \; \mathrm{and} \;  \mathbf{b}(\mathbf{x}) = \mathbf{0} 
\end{align}
\end{enumerate}

The flow of fluids under the first two classes of boundary value problems is pressure-driven, while the flow is velocity-driven for the last two. Even a few other types of boundary value problems, which in their original form do not fall into any of these classes, can be recast into one of the above classes. For example, consider a boundary value problem with a non-zero but  constant pressure applied on a part of the boundary; one can transform this problem---by shifting the datum of the pressure to the prescribed (constant) pressure---so that it belongs to Class D. 

As shown in the next section, design sensitivities can be calculated analytically for the four classes of boundary value problems.

%**********************************************;
%                                              ;
%  NAME                                        ;
%    S4_Posing_Adjoint_method.tex              ;
%                                              ;
%  WRITTEN BY                                  ;
%    Kalyana B. Nakshatrala                    ;
%                                              ;
%**********************************************;
\section{CALCULATING SENSITIVITIES USING THE ADJOINT STATE METHOD}
\label{Sec:S4_Posing_Adjoint_method}
We now estimate the sensitivity of the total dissipation rate to the permeability field. Mathematically, the task at hand is to find $D \Phi_{\mathrm{D}}[k(\mathbf{x})]$. For this estimate, we use the adjoint state method, which offers a systematic procedure to calculate sensitivities of an objective function to a set of design variables. 

The objective function, in general, depends explicitly on a set of design parameters and implicitly on the solution fields, which satisfy the forward problem. In addition, the definition of the forward problem involves the design parameters. Thus, the sensitivity of the objective function depends on the sensitivity of the solution fields. The adjoint state method circumvents an explicit calculation of the said sensitivity of the solution fields by introducing Lagrange multipliers that are the solutions of an adjoint problem---a system of partial differential equations. 

To crystallize the concept outlined above, we will relate the above description to our problem. As mentioned earlier, our task is to find the sensitivity of the total dissipation rate to the permeability field. The dissipation rate is a function of the permeability field and the velocity of the fluid. The fluid's velocity is a solution of, for example, Darcy equations---a  system of partial differential equations. Darcy equations, in turn, are written in terms of permeability. Thus, the objective function is the total dissipation rate, the design parameter is the permeability, and the state equations are Darcy equations. In the rest of this section, we derive the associated adjoint problem and calculate the desired sensitivities.  

%============================================;
%  Subsection: Deriving the adjoint problem  ;
%============================================;
\subsection{Deriving the adjoint problem} The sensitivity of the total dissipation rate functional with respect to its argument---the permeability field---takes the following mathematical form:
%---------------------------------------;
%  Equation: Definition of sensitivity  ;
%---------------------------------------;
\begin{align}
D \Phi_{\mathrm{D}}[k(\mathbf{x})] 
&= \int_{\Omega} \varphi^{\prime}_{\mathrm{D}} \Big(\mathbf{x},k(\mathbf{x}),\mathbf{v}_{\mathrm{D}}\big(\mathbf{x};k(\mathbf{x})\big)\Big) \, \mathrm{d} \Omega 
= \int_{\Omega} \left(\frac{\mu}{k(\mathbf{x})} \mathbf{v}_{\mathrm{D}}\big(\mathbf{x};k(\mathbf{x})\big) \cdot \mathbf{v}_{\mathrm{D}}\big(\mathbf{x};k(\mathbf{x})\big) \right)^{\prime} \, \mathrm{d} \Omega 
\end{align}
Following the adjoint state method, we augment the above integral with additional terms, which are essentially zero. These new terms are in the form of integrals constructed from the derivative of the state equations with respect to the design variable (i.e., $k(\mathbf{x})$). After augmenting these terms, the sensitivity can be written as follows: 
%---------------------------------------------------------;
%  Equation: Sensitivities based on adjoint state method  ;
%---------------------------------------------------------;
\begin{align}
D \Phi_{\mathrm{D}}[k(\mathbf{x})] 
&= \int_{\Omega} \left(\frac{\mu}{k(\mathbf{x})} \mathbf{v}_{\mathrm{D}}\big(\mathbf{x};k(\mathbf{x})\big) \cdot \mathbf{v}_{\mathrm{D}}\big(\mathbf{x};k(\mathbf{x})\big) \right)^{\prime} \, \mathrm{d} \Omega \notag \\
& \quad \quad \quad -2 \int_{\Omega} \boldsymbol{\Lambda}(\mathbf{x}) \cdot \Big\{\underbrace{\frac{\mu}{k(\mathbf{x})} \mathbf{v}_{\mathrm{D}}\big(\mathbf{x};k(\mathbf{x})\big) 
+ \mathrm{grad}\big[p_{\mathrm{D}}\big(\mathbf{x};k(\mathbf{x})\big)\big] - \rho \mathbf{b}(\mathbf{x})}_{Eq.~\eqref{Eqn:Posing_Darcy_BoLM}}\Big\}^{\prime} \mathrm{d} \Omega \notag \\
&\quad \quad \quad +2 \int_{\Omega} \lambda(\mathbf{x}) \,  \Big\{\underbrace{\mathrm{div}\big[\mathbf{v}_{\mathrm{D}}\big(\mathbf{x};k(\mathbf{x})\big)\big]}_{Eq.~\eqref{Eqn:Posing_Darcy_Continuity}}\Big\}^{\prime} \, \mathrm{d} \Omega \notag \\
&\quad \quad \quad +2 \int_{\Gamma^{p}} \boldsymbol{\Lambda}(\mathbf{x}) \cdot \widehat{\mathbf{n}}(\mathbf{x}) \Big\{\underbrace{p_{\mathrm{D}}\big(\mathbf{x};k(\mathbf{x})\big) - p^{\mathrm{p}}(\mathbf{x})}_{Eq.~\eqref{Eqn:Posing_Darcy_pBC}}\Big\}^{\prime} \mathrm{d} \Gamma \notag \\
&\quad \quad \quad -2 \int_{\Gamma^{v}} \lambda(\mathbf{x}) \Big\{\underbrace{\mathbf{v}_{\mathrm{D}}\big(\mathbf{x};k(\mathbf{x})\big) \cdot \widehat{\mathbf{n}}(\mathbf{x}) - v_n^{\mathrm{p}}(\mathbf{x})}_{Eq.~\eqref{Eqn:Posing_Darcy_vBC}}\Big\}^{\prime} \mathrm{d} \Gamma 
\end{align}
where $\boldsymbol{\Lambda}(\mathbf{x})$ and $\lambda(\mathbf{x})$ are the newly introduced adjoint variables. The terms in the curly brackets are the residuals of the state equations. The factor $2$ is introduced in the last four integrals for convenience, as it simplifies the adjoint problem and the expressions for the associated solution. Noting that the coefficient of viscosity $\mu$, the outward unit normal vector $\widehat{\mathbf{n}}(\mathbf{x})$, and the prescribed fields---$\rho \mathbf{b}(\mathbf{x})$, $v_{n}^{\mathrm{p}}(\mathbf{x})$ and $p^{\mathrm{p}}(\mathbf{x})$---are all independent of the permeability field, we get the following: 
%-------------------------------;
%  Equation: Simplification #1  ;
%-------------------------------;
\begin{align}
D \Phi_{\mathrm{D}}[k(\mathbf{x})] 
&= -\int_{\Omega} \frac{\mu}{k^{2}(\mathbf{x})} \mathbf{v}_{\mathrm{D}}\big(\mathbf{x};k(\mathbf{x})\big) \cdot \mathbf{v}_{\mathrm{D}}\big(\mathbf{x};k(\mathbf{x})\big) \, \mathrm{d} \Omega 
+ 2 \int_{\Omega} \frac{\mu}{k(\mathbf{x})} \mathbf{v}_{\mathrm{D}}\big(\mathbf{x};k(\mathbf{x})\big) \cdot \mathbf{v}^{\prime}_{\mathrm{D}}\big(\mathbf{x};k(\mathbf{x})\big)  \, \mathrm{d} \Omega
\notag \\
& \quad \quad \quad -2 \int_{\Omega} \boldsymbol{\Lambda}(\mathbf{x}) \cdot \Big(-\frac{\mu}{k^{2}(\mathbf{x})} \mathbf{v}_{\mathrm{D}}\big(\mathbf{x};k(\mathbf{x})\big)
+ \frac{\mu}{k(\mathbf{x})} \mathbf{v}^{\prime}_{\mathrm{D}}\big(\mathbf{x};k(\mathbf{x})\big)
+ \mathrm{grad}\big[p^{\prime}_{\mathrm{D}}\big(\mathbf{x};k(\mathbf{x})\big)\big] \Big) \mathrm{d} \Omega \notag \\
&\quad \quad \quad +2 \int_{\Omega} \lambda(\mathbf{x}) \, \mathrm{div}\big[\mathbf{v}^{\prime}_{\mathrm{D}}\big(\mathbf{x};k(\mathbf{x})\big)\big] \, \mathrm{d} \Omega \notag \\
&\quad \quad \quad +2 \int_{\Gamma^{p}} \boldsymbol{\Lambda}(\mathbf{x}) \cdot \widehat{\mathbf{n}}(\mathbf{x}) \, p^{\prime}_{\mathrm{D}}\big(\mathbf{x};k(\mathbf{x})\big) \, \mathrm{d} \Gamma \notag \\
&\quad \quad \quad 
-2 \int_{\Gamma^{v}} \lambda(\mathbf{x}) \mathbf{v}^{\prime}_{\mathrm{D}}\big(\mathbf{x};k(\mathbf{x})\big) \cdot \widehat{\mathbf{n}}(\mathbf{x}) \, \mathrm{d} \Gamma 
\end{align}
By using the Green's identity on the third term of the third integral and on the fourth integral, and grouping the terms we get the following: 
%--------------------------------;
%  Equation: Grouping the terms  ;
%--------------------------------;
\begin{align}
D \Phi_{\mathrm{D}}[k(\mathbf{x})] 
&= \int_{\Omega} \frac{\mu}{k^2(\mathbf{x})} \mathbf{v}_{\mathrm{D}}\big(\mathbf{x};k(\mathbf{x})\big) \cdot \Big(2 \boldsymbol{\Lambda}(\mathbf{x}) - \mathbf{v}_{\mathrm{D}}\big(\mathbf{x};k(\mathbf{x})\big) \Big) \, \mathrm{d} \Omega \notag \\
&\quad \quad \quad - 2 \int_{\Omega} \mathbf{v}_{\mathrm{D}}^{\prime}\big(\mathbf{x};k(\mathbf{x})\big) \cdot \left(\frac{\mu}{k(\mathbf{x})} \boldsymbol{\Lambda}(\mathbf{x}) 
    + \mathrm{grad}[\lambda(\mathbf{x})] - \frac{\mu}{k(\mathbf{x})} \mathbf{v}_{\mathrm{D}}\big(\mathbf{x};k(\mathbf{x})\big) \right) \mathrm{d} \Omega \notag \\
&\quad \quad \quad +2 \int_{\Omega} p_{\mathrm{D}}^{\prime}\big(\mathbf{x};k(\mathbf{x})\big) \,  \mathrm{div}[\boldsymbol{\Lambda}(\mathbf{x})] \, \mathrm{d} \Omega 
-2 \int_{\Gamma^{v}} 
p_{\mathrm{D}}^{\prime}\big(\mathbf{x};k(\mathbf{x})\big) \,  \mathbf{\Lambda}(\mathbf{x}) \cdot 
\widehat{\mathbf{n}}(\mathbf{x}) \, \mathrm{d} \Gamma \notag \\
&\quad \quad \quad +2 \int_{\Gamma^{p}} \mathbf{v}_{\mathrm{D}}^{\prime}\big(\mathbf{x};k(\mathbf{x})\big) \cdot \widehat{\mathbf{n}}(\mathbf{x}) \, 
\lambda(\mathbf{x}) \, \mathrm{d} \Gamma 
\end{align}

In the above equation, we do \emph{not} know \emph{a priori} the sensitivities of the solution field variables (i.e., $\mathbf{v}^{\prime}_{\mathrm{D}}(\mathbf{x};k(\mathbf{x}))$ and $p^{\prime}_{\mathrm{D}}(\mathbf{x};k(\mathbf{x}))$).
Until now the choice for the adjoint field variables is arbitrary; however, a judicious selection will circumvent the need to find the sensitivities of the solution field, and simplify the expression for design sensitivities. Specifically, the last four integral terms in the above equation can be eliminated if the adjoint field variables satisfy the following boundary value problem: 
%-----------------------------------;
%  Equation: Darcy adjoint problem  ;
%-----------------------------------;
\begin{subequations}
\label{Eqn:Posing_Adjoint_problem}
\begin{alignat}{2}
    \label{Eqn:Posing_Adjoint_problem_BoLM}
    &\frac{\mu}{k(\mathbf{x})} \boldsymbol{\Lambda}(\mathbf{x}) 
    + \mathrm{grad}[\lambda(\mathbf{x})] = \frac{\mu}{k(\mathbf{x})} \mathbf{v}_{\mathrm{D}}\big(\mathbf{x};k(\mathbf{x})\big) 
    && \quad \mathrm{in} \; \Omega \\
    \label{Eqn:Posing_Adjoint_problem_Continuity}
    &\mathrm{div}[\boldsymbol{\Lambda}(\mathbf{x})] = 0 
    && \quad \mathrm{in} \; \Omega \\
    \label{Eqn:Posing_Adjoint_problem_pBC}
    &\lambda(\mathbf{x}) = 0 
    && \quad \mathrm{on} \; \Gamma^{p} \\
    \label{Eqn:Posing_Adjoint_problem_vBC}
    &\boldsymbol{\Lambda}(\mathbf{x}) \cdot \widehat{\mathbf{n}}(\mathbf{x}) = 0 
    && \quad \mathrm{on} \; \Gamma^{v} 
\end{alignat}
\end{subequations}
Under this specific choices for the adjoint variables, the sensitivity of the objective function with respect to the permeability field reduces to:
%----------------------------------------;
%  Equation: Final design sensitivities  ;
%----------------------------------------;
\begin{align}
\label{Eqn:Posing_Darcy_final_sensitivity}
D \Phi_{\mathrm{D}}[k(\mathbf{x})] &= \int_{\Omega} \frac{\mu}{k^2(\mathbf{x})} 
\mathbf{v}_{\mathrm{D}}\big(\mathbf{x};k(\mathbf{x})\big) \cdot 
\Big(2 \boldsymbol{\Lambda}(\mathbf{x}) - \mathbf{v}_{\mathrm{D}}\big(\mathbf{x};k(\mathbf{x})\big)\Big) \, \mathrm{d} \Omega 
\end{align}

The above boundary value problem, in terms of the adjoint variables, is commonly referred as the \emph{adjoint problem}. The above set of equations governing the adjoint variables is similar to Darcy equations with a pseudo body force $\mu \mathbf{v}_{\mathrm{D}}(\mathbf{x};k(\mathbf{x}))/k(\mathbf{x})$ and homogeneous boundary conditions. In general, one has to use a numerical method to solve the adjoint problem, similar to the case with the forward problem. However, as shown below, for the four classes of boundary value problems defined earlier, we can find analytical solutions for the adjoint variables, and hence express the design sensitivities in terms of the solution fields.  

%=========================;
%  Subsection: Solutions  ;
%=========================;
\subsection{Solutions to the adjoint problem and sensitivities for the four classes} Our task is to find the adjoint variables---$\boldsymbol{\Lambda}(\mathbf{x})$ and $\lambda(\mathbf{x})$---that satisfy equations \eqref{Eqn:Posing_Adjoint_problem_BoLM}--\eqref{Eqn:Posing_Adjoint_problem_vBC}. 
%=======================================;
%  Subsubsection: Classes A and B BVPs  ;
%---------------------------------------;
\subsubsection{Classes A and B boundary value problems}
A closer inspection of Eq.~\eqref{Eqn:Posing_Adjoint_problem_BoLM} suggests that 
%---------------------------------------------;
%  Equation: Solution of the adjoint problem  ;
%---------------------------------------------;
\begin{align}
\label{Eqn:Posing_Adjoint_Problem_Classes_A_and_B_solution}
    \boldsymbol{\Lambda}(\mathbf{x}) = \mathbf{v}_{\mathrm{D}}\big(\mathbf{x};k(\mathbf{x})\big) 
    \quad \mathrm{and} \quad 
    \lambda(\mathbf{x}) = 0
\end{align}
is a plausible solution for the adjoint problem. Clearly, the above fields satisfy Eq.~\eqref{Eqn:Posing_Adjoint_problem_BoLM}. $\boldsymbol{\Lambda}(\mathbf{x})$ is divergence-free, satisfying Eq.~\eqref{Eqn:Posing_Adjoint_problem_Continuity}, as the Darcy velocity is divergence free in keeping with the continuity equation of the forward problem \eqref{Eqn:Posing_Darcy_Continuity}. For Class A problems, $\Gamma^{p} = \emptyset$, and hence Eq.~\eqref{Eqn:Posing_Adjoint_problem_vBC} is trivially satisfied. For Class B problems, $v^{\mathrm{p}}_{n}(\mathbf{x}) = 0$ on $\Gamma^{v}$, implying that 
\[
\boldsymbol{\Lambda}(\mathbf{x}) \cdot \widehat{\mathbf{n}}(\mathbf{x}) = \mathbf{v}_{\mathrm{D}}(\mathbf{x};k(\mathbf{x})) \cdot \widehat{\mathbf{n}}(\mathbf{x}) = v_n^{\mathrm{p}}(\mathbf{x}) = 0 
\]
thereby satisfying Eq.~\eqref{Eqn:Posing_Adjoint_problem_vBC}. Finally, Eq.~\eqref{Eqn:Posing_Adjoint_problem_pBC} is trivially satisfied for both the classes of problems. Thus, the above fields \eqref{Eqn:Posing_Adjoint_Problem_Classes_A_and_B_solution} offer a solution of the adjoint problem. Appealing to uniqueness of solutions for Darcy-equations, these fields are the \emph{only} solution of the adjoint problem. 

Accordingly, the sensitivity of the total dissipation rate to the permeability field, given by Eq.~\eqref{Eqn:Posing_Darcy_final_sensitivity}, for the Classes A and B problems can be simplified as follows: 
\begin{align}
    D\Phi_{\mathrm{D}}[k(\mathbf{x})] = \int_{\Omega} 
    \frac{\mu}{k^{2}(\mathbf{x})} 
    \mathbf{v}_{\mathrm{D}}(\mathbf{x};k(\mathbf{x})) \cdot \mathbf{v}_{\mathrm{D}}(\mathbf{x};k(\mathbf{x})) \, \mathrm{d} \Omega 
\end{align}
Since $\mu >0$, we have 
\begin{align}
    D\Phi_{\mathrm{D}}[k(\mathbf{x})] > 0
\end{align}
for any nonzero Darcy velocity field.

%===============================;
%  Subsubsection: Class C BVPs  ;
%-------------------------------;
\subsubsection{Class C boundary value problems}

The solution for the adjoint variables will be:
\begin{align}
\label{Eqn:Posing_solution_Class_C}
    \boldsymbol{\Lambda}(\mathbf{x}) = \mathbf{0} 
    \quad \mathrm{and} \quad 
    \lambda(\mathbf{x}) = - p_{\mathrm{D}}(\mathbf{x};k(\mathbf{x})) - \psi(\mathbf{x}) + \mathrm{constant} 
\end{align}
where $\psi(\mathbf{x})$ is the scalar potential for the body force. Recall that the pressure field can be determined up to an arbitrary constant when velocity boundary conditions are prescribed on the entire boundary---which is the case with Class C boundary value problems. The above fields trivially satisfy Eqs.~\eqref{Eqn:Posing_Adjoint_problem_Continuity} and \eqref{Eqn:Posing_Adjoint_problem_vBC}. Since $\Gamma^{p} = \emptyset$, Eq.~\eqref{Eqn:Posing_Adjoint_problem_pBC} is also met. Upon substituting the above fields in  Eq.~\eqref{Eqn:Posing_Adjoint_problem_BoLM}, we get: 
\begin{align}
    -\mathrm{grad}[p_{\mathrm{D}}(\mathbf{x};k(\mathbf{x}))] - \mathrm{grad}[\psi(\mathbf{x})] = \frac{\mu}{k(\mathbf{x})} \mathbf{v}_{\mathrm{D}}(\mathbf{x};k(\mathbf{x}))
\end{align}
Noting that $\rho \mathbf{b}(\mathbf{x}) = -\mathrm{grad}[\psi(\mathbf{x})]$, the above equation can be rewritten as:
\begin{align}
\frac{\mu}{k(\mathbf{x})} \mathbf{v}_{\mathrm{D}}(\mathbf{x};k(\mathbf{x}))
    + \mathrm{grad}[p_{\mathrm{D}}(\mathbf{x};k(\mathbf{x}))] = \rho \mathbf{b}(\mathbf{x}) 
\end{align}
which is equivalent to the balance of linear momentum equation \eqref{Eqn:Posing_Darcy_BoLM} under the forward problem. By definition, the Darcy velocity and the associated pressure satisfy Eq.~\eqref{Eqn:Posing_Darcy_BoLM}. Hence, fields given in Eq.~\eqref{Eqn:Posing_solution_Class_C} are \emph{the} solution of the adjoint problem for  Class C boundary value problems, appealing again to the uniqueness theorem for Darcy equations. 

Thus, for Class C boundary value problems, we have  
\begin{align}
    D\Phi_{\mathrm{D}}[k(\mathbf{x})] = - \int_{\Omega} 
    \frac{\mu}{k^{2}(\mathbf{x})} \mathbf{v}_{\mathrm{D}}(\mathbf{x};k(\mathbf{x})) \cdot \mathbf{v}_{\mathrm{D}}(\mathbf{x};k(\mathbf{x})) \, \mathrm{d} \Omega 
\end{align}
The positivity of the coefficient of viscosity implies that 
\begin{align}
    D\Phi_{\mathrm{D}}[k(\mathbf{x})] \leq 0
\end{align}
The above inequality becomes strictly negative for any nonzero Darcy velocity field. 

%============================;
%  Subsection: Class D BVPs  ;
%----------------------------;
\subsubsection{Class D boundary value problems}
The solution for the adjoint variables will be:
\begin{alignat}{2}
\label{Eqn:Posing_solution_Class_D}
    \boldsymbol{\Lambda}(\mathbf{x}) = \mathbf{0} 
    \quad \mathrm{and} \quad 
    \lambda(\mathbf{x}) = - p_{\mathrm{D}}(\mathbf{x};k(\mathbf{x})) 
\end{alignat}
Clearly, $\boldsymbol{\Lambda}(\mathbf{x}) = \boldsymbol{0}$ satisfies Eqs.~\eqref{Eqn:Posing_Adjoint_problem_Continuity} and \eqref{Eqn:Posing_Adjoint_problem_vBC}. Since for this class of problems, 
$p_{\mathrm{D}}(\mathbf{x};k(\mathbf{x})) = p^{p}(\mathbf{x}) = 0$ on $\Gamma^{p}$, the above solution satisfies Eq.~\eqref{Eqn:Posing_Adjoint_problem_pBC}. The above fields make the left side of Eq.~\eqref{Eqn:Posing_Adjoint_problem_DB_BoLM} to be $-\mathrm{grad}[p_{\mathrm{D}}(\mathbf{x};k(\mathbf{x}))]$, which is equal to 
\[
\frac{\mu}{k(\mathbf{x})} \mathbf{v}_{\mathrm{D}}(\mathbf{x};k(\mathbf{x}))
\] 
using  Eq.~\eqref{Eqn:Posing_Darcy_BoLM} and noting that the body force is zero for this class of problems, satisfying Eq.~\eqref{Eqn:Posing_Adjoint_problem_BoLM}. Thus, the fields in Eq.~\eqref{Eqn:Posing_solution_Class_D} are a solution of the adjoint problem; by appealing to the uniqueness theorem, these fields are the only solution for Class D boundary value problems.  

Similar to the previous class,  the sensitivity of the total dissipation rate to the permeability is: 
\begin{align}
    D\Phi_{\mathrm{D}}[k(\mathbf{x})] = - \int_{\Omega} 
    \frac{\mu}{k^{2}(\mathbf{x})} \mathbf{v}_{\mathrm{D}}\big(\mathbf{x};k(\mathbf{x})\big) \cdot \mathbf{v}_{\mathrm{D}}\big(\mathbf{x};k(\mathbf{x})\big) \, \mathrm{d} \Omega 
\end{align}
Since $\mu >0$, for any nonzero Darcy velocity field we have 
\begin{align}
    D\Phi_{\mathrm{D}}[k(\mathbf{x})] < 0
\end{align}

%==========================;
%  Subsection: Discussion  ;
%==========================;
\subsection{Discussion} For the boundary value problems under Classes A and B, the sensitivity of the total dissipation rate to the permeability field is positive. On the other hand, for boundary value problems under Classes C and D, the reverse trend is true. Mathematically, 
%------------------------------------------------------;
%  Equation: Summarizing adjoint state method results  ;
%------------------------------------------------------;
\begin{subequations} 
\begin{align}
    &\mbox{Classes A and B:} \quad 
    D\Phi_{\mathrm{D}}[k(\mathbf{x})] > 0 \\
    &\mbox{Classes C and D:} \quad 
    D\Phi_{\mathrm{D}}[k(\mathbf{x})] < 0 
\end{align}
\end{subequations}

Said differently, what we have established thus far is the total dissipation rate increases with an increase in permeability for the boundary value problems under Classes A and B; in contrast, the reverse trend holds for boundary value problems under Classes C and D. 

%**********************************************;
%                                              ;
%  NAME                                        ;
%    S5_Posing_Brinkman.tex                    ;
%                                              ;
%  WRITTEN BY                                  ;
%    Kalyana B. Nakshatrala                    ;
%                                              ;
%**********************************************;
\section{EXTENSION TO DARCY-BRINKMAN EQUATIONS}
\label{Sec:S5_Posing_Brinkman}
The Darcy-Brinkman model---a modification to the Darcy model proposed by \cite{brinkman1949calculation}---accounts for internal dissipation within the fluid besides the drag at the liquid and solid interface, considered in the Darcy model. We now show that the results established for the Darcy model extend to the Darcy-Brinkman model.

We denote the  symmetric part of gradient of a vector field $\mathbf{w}(\mathbf{x})$ as follows: 
%----------------------------------------;
%  Equation: Symmetric part of gradient  ;
%----------------------------------------;
\begin{align}
    \mathbf{D}[\mathbf{w}(\mathbf{x})] := \frac{1}{2}\left( \mathrm{grad}[\mathbf{w}(\mathbf{x})] + \mathrm{grad}[\mathbf{w}(\mathbf{x})]^{\mathrm{T}} \right)
\end{align}
where the superscript `T' denotes the transpose of a second-order tensor. We also define the following projection tensor: 
%-------------------------------;
%  Equation: Projection tensor  ;
%-------------------------------;
\begin{align}
    \label{Eqn:Posing_Projection_tensor}
    \mathbb{P}_{\|} = \mathbf{I} - \widehat{\mathbf{n}}(\mathbf{x}) \otimes 
    \widehat{\mathbf{n}}(\mathbf{x})
\end{align}
where $\mathbf{I}$ is the second-order tensor. It is easy to check that 
\begin{align}
\label{Eqn:Posing_Projection_tensor_properties}
    \mathbb{P}_{\|}^{\mathrm{T}} = \mathbb{P}_{\|}, \; 
    \mathbb{P}_{\|} \mathbb{P}_{\|} = \mathbb{P}_{\|} \; 
    \; \; \mathrm{and} \; \; 
    \mathbb{P}_{\|} \widehat{\mathbf{n}}(\mathbf{x}) = \mathbf{0} 
\end{align}
The second equation is the definition for a tensor to be a projection. The third equation implies that the normal vector  $\widehat{\mathbf{n}}(\mathbf{x})$ is not in the range space of the projection tensor $\mathbb{P}_{\|}$; said differently,  $\widehat{\mathbf{n}}(\mathbf{x})$ is in the null space of $\mathbb{P}_{\|}$.
Using this projection tensor, any vector field $\mathbf{w}(\mathbf{x})$ can be uniquely decomposed into:
%---------------------------;
%  Equation: Decomposition  ;
%---------------------------;
\begin{align}
\mathbf{w}(\mathbf{x}) = \Big(\mathbf{w}(\mathbf{x}) \cdot \widehat{\mathbf{n}}(\mathbf{x})\Big) \, 
\widehat{\mathbf{n}}(\mathbf{x})    
+ \mathbb{P}_{\|} \mathbf{w}(\mathbf{x}) 
\end{align}

Using the above notation, Darcy-Brinkman equations can be written as follows: 
%--------------------------------------;
%  Equation: Darcy-Brinkman equations  ;
%--------------------------------------;
\begin{subequations}
\begin{alignat}{2}
    \label{Eqn:Posing_DB_BoLM} 
    &\frac{\mu}{k(\mathbf{x})} \mathbf{v} + \mathrm{grad}[p] - \mathrm{div}\big[2\mu \mathbf{D}[\mathbf{v}]\big]= \rho \mathbf{b}(\mathbf{x}) 
    &&\quad \mathrm{in} \; \Omega \\
    \label{Eqn:Posing_DB_Continuity} 
    &\mathrm{div}[\mathbf{v}] = 0 
    &&\quad \mathrm{in} \; \Omega \\
    \label{Eqn:Posing_DB_Gammap_pBC} 
    &\widehat{\mathbf{n}}(\mathbf{x}) \cdot (-p \mathbf{I} + 2 \mu \mathbf{D}[\mathbf{v}])\widehat{\mathbf{n}}(\mathbf{x}) = -p^{\mathrm{p}}(\mathbf{x}) 
    &&\quad \mathrm{on} \; \Gamma^{p} \\
    \label{Eqn:Posing_DB_Gammap_vBC} 
    &\mathbb{P}_{\|}  \mathbf{v} =  \mathbf{v}^{\mathrm{p}}_{\|}(\mathbf{x}) 
    &&\quad \mathrm{on} \; \Gamma^{p} \\
    \label{Eqn:Posing_DB_vBC} 
    &\mathbf{v} = \mathbf{v}^{\mathrm{p}}(\mathbf{x}) 
    &&\quad \mathrm{on} \; \Gamma^{v} 
\end{alignat}
\end{subequations}
where $p^{\mathrm{p}}(\mathbf{x})$ is the prescribed pressure, $\mathbf{v}_{\|}^{\mathrm{p}}(\mathbf{x})$ is the prescribed tangential component of the velocity vector on $\Gamma^{p}$, and $\mathbf{v}^{\mathrm{p}}(\mathbf{x})$ is the prescribed (full) velocity vector on $\Gamma^{v}$. 

If velocity boundary conditions are prescribed on the entire boundary (i.e., $\Gamma^{v} = \partial \Omega$), similar to Darcy equations, solutions to Darcy-Brinkman equations exist only if the following compatibility condition is satisfied: 
%-------------------------------------;
%  Equation: Compatibility condition  ;
%-------------------------------------;
\begin{align}
    \label{Eqn:Posing_Compatibility_DB}
    \int_{\Gamma^{v} = \partial \Omega}
    \mathbf{v}^{\mathrm{p}}(\mathbf{x}) \cdot \widehat{\mathbf{n}}(\mathbf{x}) \, \mathrm{d} \Gamma = 0 
\end{align}
The solutions to Darcy-Brinkman equations are unique except when $\Gamma^{v} = \partial \Omega$, where, upon meeting the above compatibility condition \eqref{Eqn:Posing_Compatibility_DB}, the pressure is unique only up to an arbitrary constant. See \citep{shabouei_nakshatrala_cicp} for a proof of the unique theorem for Darcy-Brinkman equations. From hereon, the solutions fields---the Darcy-Brinkman velocity and the associated pressure---are denoted as $\mathbf{v}_{\mathrm{DB}}(\mathbf{x};k(\mathbf{x}))$ and 
$p_{\mathrm{DB}}(\mathbf{x};k(\mathbf{x}))$, respectively. 

The dissipation rate density for an arbitrary vector field $\mathbf{v}(\mathbf{x})$ under the Darcy-Brinkman model takes the following mathematical form: 
%--------------------------------------;
%  Equation: Dissipation rate density  ;
%--------------------------------------;
\begin{align}
\varphi_{\mathrm{DB}}\big(\mathbf{x},k(\mathbf{x}),\mathbf{v}(\mathbf{x})\big) 
= \frac{\mu}{k(\mathbf{x})} \mathbf{v} (\mathbf{x}) \cdot \mathbf{v}(\mathbf{x})
+ 2 \mu \,  \mathbf{D}[\mathbf{v}(\mathbf{x})] \cdot \mathbf{D}[\mathbf{v}(\mathbf{x})]
\end{align}
The total dissipation rate functional---to denote the total dissipation rate expended by the \emph{Darcy-Brinkman velocity} under the \emph{Darcy-Brinkman model}---is defined as: 
%------------------------------------;
%  Equation: Total dissipation rate  ;
%------------------------------------;
\begin{align}
\Phi_{\mathrm{DB}}[k(\mathbf{x})] 
:= \int_{\Omega} 
\varphi_{\mathrm{DB}}\big(\mathbf{x},k(\mathbf{x}),\mathbf{v}_{\mathrm{DB}}(\mathbf{x};k(\mathbf{x}))\big) \, \mathrm{d} \Omega 
\end{align}

%===================================;
%  Subsection: Four classes of BVP  ;
%===================================;
\subsection{Four classes of boundary value problems under the Darcy-Brinkman model} 
Since the boundary conditions under the Darcy-Brinkman model differ from that of the Darcy model, we need to slightly modify the definitions of the classes of the boundary value problems. Still, the flow of fluids under Classes A and B will be pressure driven while Classes C and D consider velocity-driven flows. 
%-------------------------------------;
%  Enumerate: Classification of BVPs  ;
%-------------------------------------;
\begin{enumerate}
\item \textbf{Class A:} Prescription of a pressure loading on the entire boundary with the corresponding tangential velocity to be zero. Stated mathematically, 
\begin{align}
\Gamma^{v} = \emptyset 
\; \mathrm{and} \; 
\mathbf{v}_{\|}(\mathbf{x}) = \mathbf{0} \; \mathrm{on} \; \Gamma^{p} = \partial \Omega 
\end{align}
\item \textbf{Class B:} Prescription of pressure loading on a part of the boundary while homogeneous velocity boundary conditions on the rest. That is, 
\begin{align}
\Gamma^{p} \neq \emptyset, \; 
\Gamma^{v} \neq \emptyset, 
\mathbf{v}_{\|}^{\mathrm{p}}(\mathbf{x}) = \mathbf{0} \; \mathrm{on} \;  \Gamma^{p}, 
\; \mathrm{and} \;
\mathbf{v}^{\mathrm{p}}(\mathbf{x}) = \mathbf{0} \; \mathrm{on} \;  \Gamma^{v}
\end{align}
\item \textbf{Class C:} Prescription of a compatible velocity boundary condition on the entire boundary and requiring the prescribed body force to be a conservative vector field. Mathematically, 
\begin{align}
\Gamma^{p} = \emptyset, \; \int_{\Gamma^{v} = \partial \Omega} \mathbf{v}^{\mathrm{p}}(\mathbf{x}) \cdot \widehat{\mathbf{n}}(\mathbf{x}) \, \mathrm{d} \Gamma = 0 \quad \mathrm{and} \quad \rho\mathbf{b}(\mathbf{x}) = - \mathrm{grad}[\psi] 
\end{align}
where $\psi(\mathbf{x})$ is a scalar field. 
\item \textbf{Class D:} Prescription of velocity boundary conditions on a portion of the boundary and homogeneous (pressure and tangential velocity) boundary conditions on the rest, and requiring zero body force. Written mathematically, 
\begin{align}
\Gamma^{p} \neq \emptyset, \; \Gamma^{v} \neq 0, \;  
p^{\mathrm{p}}(\mathbf{x}) = 0 
\; \mathrm{and} \; \mathbf{v}_{\|}^{\mathrm{p}}(\mathbf{x}) = \mathbf{0} \; \mathrm{on} \; \Gamma^{p}, \; \mathrm{and} \;  \mathbf{b}(\mathbf{x}) = \mathbf{0} 
\end{align}
\end{enumerate}

%============================================================;
%  Subsection: Sensitivities using the adjoint state method  ;
%============================================================;
\subsection{Sensitivities using the adjoint state method}
We now calculate the sensitivity of the total dissipation rate to the permeability field, $D\Phi_{\mathrm{DB}}[k(\mathbf{x})]$, which takes the mathematical form:
%---------------------------------------------------------;
%  Equation: Sensitivities based on adjoint state method  ;
%---------------------------------------------------------;
\begin{align}
D \Phi_{\mathrm{DB}}[k(\mathbf{x})] 
&= \int_{\Omega} \varphi^{\prime}_{\mathrm{DB}}\Big(\mathbf{x},k(\mathbf{x}),\mathbf{v}_{\mathrm{DB}}\big(\mathbf{x};k(\mathbf{x})\big)\Big) \mathrm{d} \Omega \notag \\ 
&= \int_{\Omega} \left(\frac{\mu}{k(\mathbf{x})} \mathbf{v}_{\mathrm{DB}}\big(\mathbf{x};k(\mathbf{x})\big) \cdot \mathbf{v}_{\mathrm{DB}}\big(\mathbf{x};k(\mathbf{x})\big) + 2 \mu \,  \mathbf{D}\big[\mathbf{v}_{\mathrm{DB}}\big(\mathbf{x};k(\mathbf{x})\big)\big] \cdot \mathbf{D}\big[\mathbf{v}_{\mathrm{DB}}\big(\mathbf{x};k(\mathbf{x})\big)\big] \right)^{\prime} \mathrm{d} \Omega
\end{align}
Following the adjoint state method, we augment the above term with the derivative of the residuals of Darcy-Brinkman equations: 
%---------------------------------------------------------;
%  Equation: Sensitivities based on adjoint state method  ;
%---------------------------------------------------------;
{\small
\begin{align}
D \Phi_{\mathrm{DB}}[k(\mathbf{x})] 
&= \int_{\Omega} \left(\frac{\mu}{k(\mathbf{x})} \mathbf{v}_{\mathrm{DB}}\big(\mathbf{x};k(\mathbf{x})\big) \cdot \mathbf{v}_{\mathrm{DB}}\big(\mathbf{x};k(\mathbf{x})\big) 
+ 2 \mu \,  \mathbf{D}[\mathbf{v}_{\mathrm{DB}}(\mathbf{x};k(\mathbf{x}))] \cdot 
\mathbf{D}[\mathbf{v}_{\mathrm{DB}}(\mathbf{x};k(\mathbf{x}))] 
\right)^{\prime} \, \mathrm{d} \Omega \notag \\
&\quad \quad -2 \int_{\Omega} \boldsymbol{\Lambda}(\mathbf{x}) \cdot \Big\{\underbrace{\frac{\mu}{k(\mathbf{x})} \mathbf{v}_{\mathrm{DB}}(\mathbf{x};k(\mathbf{x})) 
+ \mathrm{grad}[p_{\mathrm{DB}}(\mathbf{x};k(\mathbf{x}))] - \mathrm{div}\big[2 \mu \,  \mathbf{D}[\mathbf{v}_{\mathrm{DB}}(\mathbf{x};k(\mathbf{x}))]\big]- \rho \mathbf{b}(\mathbf{x})}_{Eq.~\eqref{Eqn:Posing_DB_BoLM}}\Big\}^{\prime} \mathrm{d} \Omega \notag \\
&\quad \quad +2 \int_{\Omega} 
\lambda(\mathbf{x}) \,  
\Big\{\underbrace{\mathrm{div}[\mathbf{v}_{\mathrm{DB}}(\mathbf{x};k(\mathbf{x}))]}_{Eq.~\eqref{Eqn:Posing_DB_Continuity}}\Big\}^{\prime} \, \mathrm{d} \Omega \notag \\
&\quad \quad -2 \int_{\Gamma^{p}} \boldsymbol{\Lambda}(\mathbf{x}) \cdot \widehat{\mathbf{n}}(\mathbf{x}) \,  \Big\{\underbrace{\widehat{\mathbf{n}}(\mathbf{x})\cdot \Big(-p_{\mathrm{DB}}(\mathbf{x};k(\mathbf{x})) \mathbf{I} + 2 \mu \mathbf{D}\big[\mathbf{v}_{\mathrm{DB}}(\mathbf{x};k(\mathbf{x}))\big]\Big)\widehat{\mathbf{n}}(\mathbf{x}) + p^{\mathrm{p}}(\mathbf{x})}_{Eq.~\eqref{Eqn:Posing_DB_Gammap_pBC}}\Big\}^{\prime} 
\mathrm{d} \Gamma \notag \\ 
&\quad \quad -2 \int_{\Gamma^{p}} \left\{\mathbb{P}_{\|} \Big(2 \mu \mathbf{D}[ \boldsymbol{\Lambda}(\mathbf{x})] 
- 2 \mu \mathbf{D}[ \mathbf{v}_{\mathrm{DB}}(\mathbf{x};k(\mathbf{x}))] 
\Big) 
\widehat{\mathbf{n}}(\mathbf{x}) \right\} \cdot  \Big\{\underbrace{\mathbb{P}_{\|} \mathbf{v}_{\mathrm{DB}}(\mathbf{x};k(\mathbf{x})) -  \mathbf{v}_{\|}^{\mathrm{p}}(\mathbf{x})}_{Eq.~\eqref{Eqn:Posing_DB_Gammap_vBC}}\Big\}^{\prime} 
\mathrm{d} \Gamma \notag \\ 
&\quad \quad +2 \int_{\Gamma^{v}} \Big(-\lambda(\mathbf{x}) \mathbf{I} + 2 \mu \mathbf{D}[\boldsymbol{\Lambda}(\mathbf{x})]
- 2 \mu \mathbf{D}[\mathbf{v}_{\mathrm{DB}}(\mathbf{x};k(\mathbf{x}))]\Big)  \widehat{\mathbf{n}}(\mathbf{x}) \cdot \Big\{\underbrace{\mathbf{v}_{\mathrm{DB}}(\mathbf{x};k(\mathbf{x})) - \mathbf{v}^{\mathrm{p}}(\mathbf{x})}_{Eq.~\eqref{Eqn:Posing_DB_vBC}}\Big\}^{\prime} \mathrm{d} \Gamma
\end{align}}
The rationale behind the incorporation of the factor 2 in the last four terms and the selection of the weights within the integrals, especially in the last two terms, is to simplify the resulting adjoint problem. We note that the coefficient of viscosity $\mu$, the outward unit normal vector $\widehat{\mathbf{n}}(\mathbf{x})$, and the prescribed fields---$\rho \mathbf{b}(\mathbf{x})$, $\mathbf{v}_{\|}^{\mathrm{p}}(\mathbf{x})$,  $\mathbf{v}^{\mathrm{p}}(\mathbf{x})$ and $p^{\mathrm{p}}(\mathbf{x})$---are all independent of the permeability field. Using Green's identity and the properties of projection tensor given in Eq.~\eqref{Eqn:Posing_Projection_tensor_properties}, performing algebraic simplifications, and grouping the terms, we get the following:
%----------------------------;
%  Equation: Simplification  ;
%----------------------------;
{\small
\begin{align}
D \Phi_{\mathrm{DB}}[k(\mathbf{x})] 
&= \int_{\Omega} \frac{\mu}{k^2(\mathbf{x})} \mathbf{v}_{\mathrm{DB}}(\mathbf{x};k(\mathbf{x})) \cdot \big(2 \boldsymbol{\Lambda}(\mathbf{x}) - \mathbf{v}_{\mathrm{DB}}(\mathbf{x};k(\mathbf{x}))\big) \, \mathrm{d} \Omega \notag \\
&- 2 \int_{\Omega} \mathbf{v}_{\mathrm{DB}}^{\prime}(\mathbf{x};k(\mathbf{x})) \cdot \Big\{ \frac{\mu}{k(\mathbf{x})} \boldsymbol{\Lambda}(\mathbf{x})
+ \mathrm{grad}[\lambda(\mathbf{x})]
- \mathrm{div}\Big[2 \mu \, \mathbf{D}\big[\boldsymbol{\Lambda}(\mathbf{x})\big]\Big] \notag \\
&\qquad \qquad \qquad \qquad \qquad \qquad \qquad - \Big(\frac{\mu}{k(\mathbf{x})} \mathbf{v}_{\mathrm{DB}}(\mathbf{x};k(\mathbf{x})) - \mathrm{div}\big[2 \mu \mathbf{D}[\mathbf{v}_{\mathrm{DB}}(\mathbf{x};k(\mathbf{x}))]\big] \Big)\Big\} \, \mathrm{d} \Omega \notag \\ 
&+2 \int_{\Omega}
p_{\mathrm{DB}}^{\prime}(\mathbf{x};k(\mathbf{x})) \,  \mathrm{div}[\boldsymbol{\Lambda}(\mathbf{x})] \, \mathrm{d} \Omega \notag \\
&-2 \int_{\Gamma^{p}} \Big(
\mathbf{v}_{\mathrm{DB}}^{\prime}(\mathbf{x};k(\mathbf{x})) \cdot \widehat{\mathbf{n}}(\mathbf{x}) \Big) \,  \widehat{\mathbf{n}}(\mathbf{x}) \cdot  \Big(-\lambda(\mathbf{x}) \mathbf{I} + 
2 \mu \, \mathbf{D}[\boldsymbol{\Lambda}(\mathbf{x})] - 2 \mu \, \mathbf{D}[\mathbf{v}_{\mathrm{DB}}(\mathbf{x};k(\mathbf{x}))]\Big) \widehat{\mathbf{n}}(\mathbf{x}) \,  \mathrm{d} \Gamma \notag \\
&+2 \int_{\Gamma^{p}}
\Big(\mathbb{P}_{\|}
2 \mu \, \mathbf{D}[\mathbf{v}_{\mathrm{DB}}^{\prime}(\mathbf{x};k(\mathbf{x}))]  \widehat{\mathbf{n}}(\mathbf{x}) \Big) \cdot \Big(\mathbb{P}_{\|} \boldsymbol{\Lambda}(\mathbf{x}) \Big) \,  \mathrm{d} \Gamma \notag \\
&+2 \int_{\Gamma^{v}}
\left(-p_{\mathrm{DB}}^{\prime}(\mathbf{x};k(\mathbf{x})) \mathbf{I} + 2 \mu \mathbf{D}\big[\mathbf{v}_{\mathrm{DB}}^{\prime}(\mathbf{x};k(\mathbf{x}))\big]\right)\widehat{\mathbf{n}}(\mathbf{x}) \cdot \boldsymbol{\Lambda}(\mathbf{x}) \, \mathrm{d} \Gamma 
\end{align}}

Following the same approach as in the case of the Darcy model, we choose the adjoint variables to make the last five integrals vanish. Accordingly, $\boldsymbol{\Lambda}(\mathbf{x})$ and $\lambda(\mathbf{x})$ satisfy the following adjoint problem, corresponding to Darcy-Brinkman equations: 
%------------------------------------;
%  Equation: Adjoint problem for DB  ;
%------------------------------------;
\begin{subequations}
\label{Eqn:Posing_Adjoint_problem_DB}
\begin{alignat}{2}
    \label{Eqn:Posing_Adjoint_problem_DB_BoLM}
    &\frac{\mu}{k(\mathbf{x})} \boldsymbol{\Lambda}(\mathbf{x}) 
    + \mathrm{grad}[\lambda(\mathbf{x})] 
    - \mathrm{div}\big[2\mu   \mathbf{D}[\boldsymbol{\Lambda}(\mathbf{x})]\big] = \frac{\mu}{k(\mathbf{x})} \mathbf{v}_{\mathrm{DB}}\big(\mathbf{x};k(\mathbf{x})\big) \notag \\
    & \qquad \qquad \qquad \qquad \qquad \qquad \qquad \qquad \qquad \qquad - \mathrm{div}\big[2 \mu  \mathbf{D}[\mathbf{v}_{\mathrm{DB}}(\mathbf{x};k(\mathbf{x}))]\big]
    && \; \mathrm{in} \; \Omega \\
    &\mathrm{div}[\boldsymbol{\Lambda}(\mathbf{x})] = 0 
    && \; \mathrm{in} \; \Omega \\
    &\widehat{\mathbf{n}}(\mathbf{x}) \cdot \Big(-\lambda(\mathbf{x}) \mathbf{I} + 2 \mu \, 
    \mathbf{D}[\boldsymbol{\Lambda}(\mathbf{x})]
    \Big) \widehat{\mathbf{n}}(\mathbf{x}) = \widehat{\mathbf{n}}(\mathbf{x}) \cdot \Big(2 \mu \,  \mathbf{D}[\mathbf{v}_{\mathrm{DB}}(\mathbf{x};k(\mathbf{x}))]\Big)  \widehat{\mathbf{n}}(\mathbf{x}) 
    && \; \mathrm{on} \; \Gamma^{p} \\ 
    &\mathbb{P}_{\|}  \boldsymbol{\Lambda}(\mathbf{x}) = \mathbf{0} && \; \mathrm{on} \; \Gamma^{p} \\ 
    \label{Eqn:Posing_Adjoint_problem_DB_vBC}
    &\boldsymbol{\Lambda}(\mathbf{x}) = \mathbf{0} 
    && \; \mathrm{on} \; \Gamma^{v} 
\end{alignat}
\end{subequations}

The above adjoint problem resembles Darcy-Brinkman equations \eqref{Eqn:Posing_DB_BoLM}--\eqref{Eqn:Posing_DB_vBC} but with homogeneous velocity boundary conditions, a non-homogeneous body force of 
\[
\frac{\mu}{k(\mathbf{x})} \mathbf{v}_{\mathrm{DB}}\big(\mathbf{x};k(\mathbf{x})\big)  - \mathrm{div}\big[2 \mu  \mathbf{D}[\mathbf{v}_{\mathrm{DB}}(\mathbf{x};k(\mathbf{x}))]\big]
\]
and a pressure loading of 
\[
-\widehat{\mathbf{n}}(\mathbf{x}) \cdot \Big(2 \mu \,  \mathbf{D}[\mathbf{v}_{\mathrm{DB}}(\mathbf{x};k(\mathbf{x}))]\Big)  \widehat{\mathbf{n}}(\mathbf{x})
\]
The negative sign in the above expression indicates that a pressure loading acts inwards, opposing the unit outward normal vector $\widehat{\mathbf{n}}(\mathbf{x})$. 

For the adjoint variables satisfying the adjoint problem \eqref{Eqn:Posing_Adjoint_problem_DB_BoLM}--\eqref{Eqn:Posing_Adjoint_problem_DB_vBC}, 
the sensitivity of the total dissipation rate with respect to the permeability field is:
%----------------------------------------;
%  Equation: Final sensitivities for DB  ;
%----------------------------------------;
\begin{align}
D \Phi_{\mathrm{DB}}[k(\mathbf{x})] &= \int_{\Omega} \frac{\mu}{k^2(\mathbf{x})} 
\mathbf{v}_{\mathrm{DB}}\big(\mathbf{x};k(\mathbf{x})\big) \cdot 
\Big(2 \boldsymbol{\Lambda}(\mathbf{x}) - \mathbf{v}_{\mathrm{DB}}\big(\mathbf{x};k(\mathbf{x})\big)\Big) \, \mathrm{d} \Omega 
\end{align}
The above expression is \emph{similar} to that of the one obtained under the Darcy model (cf.~Eq.~\eqref{Eqn:Posing_Darcy_final_sensitivity}). To find the sensitivity, we need to know $\boldsymbol{\Lambda}(\mathbf{x})$, which is a solution of the adjoint problem. Note that $\mathbf{v}_{\mathrm{DB}}(\mathbf{x};k(\mathbf{x}))$ is known from the forward problem. 

%=======================================;
%  Subsubsection: Classes A and B BVPs  ;
%---------------------------------------;
\subsubsection{Classes A and B boundary value problems}
Even under the Darcy-Brinkman model, the solution for the adjoint variables for these two classes of boundary value problem will be:
\begin{align}
    \boldsymbol{\Lambda}(\mathbf{x}) = \mathbf{v}_{\mathrm{DB}}\big(\mathbf{x};k(\mathbf{x})\big) 
    \quad \mathrm{and} \quad 
    \lambda(\mathbf{x}) = 0 
\end{align}
Following the same reasoning as done for Darcy equations and appealing to the uniqueness theorem for Darcy-Brinkman equations, the above pair is the only solution of the adjoint problem \eqref{Eqn:Posing_Adjoint_problem_DB_BoLM}--\eqref{Eqn:Posing_Adjoint_problem_DB_vBC}. Consequently, 
\begin{align}
    D\Phi_{\mathrm{DB}}[k(\mathbf{x})] = \int_{\Omega} 
    \frac{\mu}{k^{2}(\mathbf{x})} 
    \mathbf{v}_{\mathrm{DB}}(\mathbf{x};k(\mathbf{x})) \cdot \mathbf{v}_{\mathrm{DB}}(\mathbf{x};k(\mathbf{x})) \, \mathrm{d} \Omega 
\end{align}
Since $\mu >0$, we have 
\begin{align}
    D\Phi_{\mathrm{DB}}[k(\mathbf{x})] > 0
\end{align}
for any nonzero Darcy-Brinkman velocity field. 

%===============================;
%  Subsubsection: Class C BVPs  ;
%-------------------------------;
\subsubsection{Class C boundary value problems}

The solution for the adjoint variables will be:
\begin{align}
    \boldsymbol{\Lambda}(\mathbf{x}) = \mathbf{0} 
    \quad \mathrm{and} \quad 
    \lambda(\mathbf{x}) = - p_{\mathrm{DB}}(\mathbf{x};k(\mathbf{x})) - \psi(\mathbf{x}) + \mathrm{constant} 
\end{align}
For this class of boundary value problems, the sensitivity of the total dissipation rate to the permeability field can be compactly written as follows:
\begin{align}
    D\Phi_{\mathrm{DB}}[k(\mathbf{x})] = - \int_{\Omega} 
    \frac{\mu}{k^{2}(\mathbf{x})} \mathbf{v}_{\mathrm{DB}}(\mathbf{x};k(\mathbf{x})) \cdot \mathbf{v}_{\mathrm{DB}}(\mathbf{x};k(\mathbf{x})) \, \mathrm{d} \Omega 
\end{align}
The coefficient of viscosity is positive, implying 
\begin{align}
    D\Phi_{\mathrm{DB}}[k(\mathbf{x})] < 0
\end{align}
for any nonzero Darcy-Brinkman velocity field. 

%===============================;
%  Subsubsection: Class D BVPs  ;
%-------------------------------;
\subsubsection{Class D boundary value problems}
The solution for the adjoint variables will be:
\begin{align}
    \boldsymbol{\Lambda}(\mathbf{x}) = \mathbf{0} 
    \quad \mathrm{and} \quad 
    \lambda(\mathbf{x}) = - p_{\mathrm{DB}}(\mathbf{x};k(\mathbf{x}))
\end{align}
The sensitivity of the total dissipation rate to the permeability field is again:
\begin{align}
    D\Phi_{\mathrm{DB}}[k(\mathbf{x})] = - \int_{\Omega} 
    \frac{\mu}{k^{2}(\mathbf{x})} \mathbf{v}_{\mathrm{DB}}\big(\mathbf{x};k(\mathbf{x})\big) \cdot \mathbf{v}_{\mathrm{DB}}\big(\mathbf{x};k(\mathbf{x})\big) \, \mathrm{d} \Omega 
\end{align}
For any nonzero Darcy-Brinkman velocity field, we have
\begin{align}
    D\Phi_{\mathrm{DB}}[k(\mathbf{x})] < 0
\end{align}
as $\mu > 0$.

%==========================;
%  Subsection: Discussion  ;
%==========================;
\subsection{Discussion} 
Although Darcy-Brinkman equations consider an additional mode of dissipation---internal friction in the fluid, the trend for the sensitivity of the total dissipation rate to the permeability field remains the same as that of Darcy equations. Moreover, expressions for the sensitivities are similar except for the apparent replacement of  $\mathbf{v}_{\mathrm{D}}(\mathbf{x};k(\mathbf{x}))$ with $\mathbf{v}_{\mathrm{DB}}(\mathbf{x};k(\mathbf{x}))$. Thus, what we have established is, for the boundary value problems under Classes A and B, the total dissipation rate increases with an increase in the permeability. On the other hand, the opposite is true for boundary value problems under Classes C and D. Mathematically, the trend can be summarized as follows: 
\begin{subequations}
\begin{align}
    &\mbox{Classes A and B:} \quad 
    k_1(\mathbf{x}) \geq k_2(\mathbf{x})  
    \quad \mathrm{implies} \quad 
    \Phi_{\mathrm{DB}}[k_{1}(\mathbf{x})] \geq \Phi_{\mathrm{DB}}[k_2(\mathbf{x})] \\
    &\mbox{Classes C and D:} \quad 
    k_1(\mathbf{x}) \geq k_2(\mathbf{x})  
    \quad \mathrm{implies} \quad 
    \Phi_{\mathrm{DB}}[k_{1}(\mathbf{x})] \leq \Phi_{\mathrm{DB}}[k_2(\mathbf{x})] 
\end{align}
\end{subequations}

%**********************************************;
%                                              ;
%  NAME                                        ;
%    S6_Posing_Material_design.tex             ;
%                                              ;
%  WRITTEN BY                                  ;
%    Kalyana B. Nakshatrala                    ;
%                                              ;
%**********************************************;
\section{POSING DESIGN PROBLEMS UNDER TOPOLOGY OPTIMIZATION}
\label{Sec:S6_Posing_Material_design}
In topology optimization, the design variable is selection of the material. In porous media applications, the material property is the permeability. The extremization can be minimization or maximization. The questions is should we minimize or maximize, and should we place a volumetric bound on the usage of the high-permeability or low-permeability material. The actual choice depends on the nature of the boundary value problem (i.e., on the forward problem). 

In a recently published paper, this question has been addressed by solving one-dimensional problems. However, no rigorous justification is provided, and we now provide the much needed theoretical justification. In general, four combinations of material design problems are possible: \emph{maximize/minimize} the total dissipation rate with a volume bound constraint on the use of \emph{high}- or \emph{low}-permeability material. 

% These sensitivities provide rigorous answers to several related design questions: (a) how to pose a material design problem---either maximize or minimize the total dissipation rate, (b) should the volumetric bound constraint be on the use of high- or low-permeability material, and (c) when does the optimization problem lead to nontrivial designs.

For Classes A and B, the total dissipation rate increases with the permeability field. For such problems, if the total dissipation rate is \emph{maximized} with a volume bound constraint only on the low-permeability material---no bound on the use of the high-permeability material---then the design is trivial: the optimal design is the one with the high-permeability material everywhere. On the other hand, for the same classes of problems, if the volume bound constraint is placed on the use of the high-permeability material, the optimal design is nontrivial. If the total dissipation rate is \emph{minimized}, then volume bound constraint on the high-permeability material produces a trivial design---the low-permeability material is placed everywhere in the domain, while a volume bound constraint is placed on the use of low-permeability material then the design is nontrivial. 

Table \ref{Table:Posing_Scenarios_BVPs} summarizes all these said scenarios---valid for both Darcy and Darcy-Brinkman models. 

%-------------------------------------------------;
%  Table: Scenarios for the four classes of BVPs  ;
%-------------------------------------------------;
\begin{table}[ht]
  \centering
  \caption{Scenarios under which trivial and nontrivial designs occur for the four classes of boundary value problems.}
    \begin{tabular}{|c|c|c|c|c|c|}
    \hline
    \multicolumn{2}{|c|}{\multirow{2}[4]{*}{}} & \multicolumn{4}{c|}{\cellcolor{yellow!25}\textbf{Extremization of the total dissipation rate}} \bigstrut\\
    \hhline{~~|*4{-}|}
    % \cline{3-6}    
    %
\multicolumn{2}{|c|}{} & \multicolumn{2}{c|}{\cellcolor{yellow!25}\emph{\textbf{Maximization}}} & \multicolumn{2}{c|}{\cellcolor{yellow!25}\emph{\textbf{Minimization}}} \bigstrut\\
    \hhline{~~|*4{-}|}
    % \cline{3-6}    
    %
    \multicolumn{2}{|c|}{} & \cellcolor{gray!20}\textbf{Classes A \& B} & \cellcolor{gray!05} \textbf{Classes C \& D} & \cellcolor{gray!20} \textbf{Classes A \& B} & \cellcolor{gray!05} \textbf{Classes C \& D}  \bigstrut\\ \hline
    %%%%
    \cellcolor{green!20} & \cellcolor{green!20} & \cellcolor{gray!20} & \cellcolor{gray!5} & \cellcolor{gray!20} & \cellcolor{gray!5} \bigstrut\\
    \cellcolor{green!20} & \cellcolor{green!20} & \cellcolor{gray!20} \textbf{nontrivial design} & \cellcolor{gray!5} \textbf{trivial design} & \cellcolor{gray!20} \textbf{trivial design} & \cellcolor{gray!5} \textbf{nontrivial design}  \bigstrut\\
    \cellcolor{green!20} & \cellcolor{green!20} & \cellcolor{gray!20} & \cellcolor{gray!5} \emph{low-permeability} & \cellcolor{gray!20} \emph{low-permeability} & \cellcolor{gray!5}  \bigstrut\\
    \cellcolor{green!20} & \cellcolor{green!20} & \cellcolor{gray!20} & \cellcolor{gray!5} \emph{material placed in} & \cellcolor{gray!20} \emph{material placed in} & \cellcolor{gray!5} \bigstrut\\
    \cellcolor{green!20} & \cellcolor{green!20} & \cellcolor{gray!20} & \cellcolor{gray!5} \emph{the entire domain} & \cellcolor{gray!20} \emph{the entire domain} & \cellcolor{gray!5}  \bigstrut\\
    \cellcolor{green!20} & \cellcolor{green!20}\multirow{-6}[10]{*}{\begin{sideways} \textbf{\emph{high-permeability}} \hspace{0.05in} 
    \end{sideways}} & \cellcolor{gray!20} & \cellcolor{gray!5} & \cellcolor{gray!20} & \cellcolor{gray!5} \bigstrut \\
    %%%%
    \hhline{~|*5{-}|}
    % \cline{2-6}  
    %%%%
    \cellcolor{green!20} & \cellcolor{green!25} & \cellcolor{gray!20} & \cellcolor{gray!5} & \cellcolor{gray!20} & \cellcolor{gray!5} \bigstrut\\
    \cellcolor{green!20} & \cellcolor{green!25} & \cellcolor{gray!20} \textbf{trivial design} & \cellcolor{gray!5} \textbf{nontrivial design} & \cellcolor{gray!20}\textbf{nontrivial design} & \cellcolor{gray!5} \textbf{trivial design}  \bigstrut\\
    \cellcolor{green!20} & \cellcolor{green!20} & \cellcolor{gray!20}\emph{high-permeability} & \cellcolor{gray!5} & \cellcolor{gray!20} & \cellcolor{gray!5} \emph{high-permeability}  \bigstrut\\
    \cellcolor{green!20} & \cellcolor{green!20} & \cellcolor{gray!20} \emph{material placed in} & \cellcolor{gray!5} & \cellcolor{gray!20} & \cellcolor{gray!5} \emph{material placed in}  \bigstrut\\
    \cellcolor{green!20} & \cellcolor{green!20} & \cellcolor{gray!20} \emph{the entire domain} & \cellcolor{gray!5} & \cellcolor{gray!20} & \cellcolor{gray!5} \emph{the entire domain}  \bigstrut\\
    \cellcolor{green!20} \multirow{-14}[10]{*}{\begin{sideways}\textbf{Volume bound constraint} \vspace{0.05in} \end{sideways}} & \cellcolor{green!20}\multirow{-6}[10]{*}{\begin{sideways} \textbf{\emph{low-permeability}} \hspace{0.05in}\end{sideways}} & \cellcolor{gray!20} &\cellcolor{gray!5} & \cellcolor{gray!20} & \cellcolor{gray!5} \bigstrut \\\hline
    %%%%
    \end{tabular}
  \label{Table:Posing_Scenarios_BVPs}
\end{table}

%**********************************************;
%                                              ;
%  NAME                                        ;
%    S7_Posing_CR.tex                          ;
%                                              ;
%  WRITTEN BY                                  ;
%    Kalyana B. Nakshatrala                    ;
%                                              ;
%**********************************************;
\section{CLOSURE}
\label{Sec:S7_Posing_CR}
This paper addressed how to pose well-posed material design problems---using the rate of dissipation---for flow through porous media applications. The answer to this design question depends on specifics of the boundary value problem---the so-called forward/direct problem. Mathematical analysis using the adjoint state method identified four classes of boundary value problems: two of them involve pressure-driven flows, and the other two are velocity-driven. Using these four classes, this paper rigorously answered:
\begin{enumerate}
\item whether to maximize or minimize the total dissipation rate, 
\item should we place a bound constraint on high- or low-permeability material, and 
\item conditions for nontrivial designs.
\end{enumerate} 
This newfound knowledge puts the material design problem using the dissipation rate under a sound footing.

While answering the design question, this study revealed a fundamental property of the flow of fluids through porous media: how the total dissipation rate varies with the permeability field for the four classes of boundary value problems. The said sensitivity is positive for the two classes with prescribed pressure loading driving the flow. However, the trend is opposite for the other two classes in which velocity boundary conditions drive the flow. Yet another use of the analytical sensitivities derived in this paper is that a computer implementation can use them without solving the adjoint problem numerically, decreasing the time to compute. These four classes of boundary value problems are fundamental, with far-reaching consequences beyond posing material design problems.

\appendix 
\setcounter{equation}{0}
\renewcommand{\theequation}{A\arabic{equation}}
%**********************************************;
%                                              ;
%  NAME                                        ;
%    Appendix.tex                              ;
%                                              ;
%  WRITTEN BY                                  ;
%    Kalyana B. Nakshatrala                    ;
%                                              ;
%**********************************************;
\section{AN ALTERNATIVE FORM OF DARCY-BRINKMAN EQUATIONS}
\label{Sec:App_Posing_Alternative_Brinkman}
An alternative form of Darcy-Brinkman equations, found in the literature, has the following mathematical form: 
%--------------------------------------;
%  Equation: Alternative DB equations  ;
%--------------------------------------;
\begin{subequations}
\begin{alignat}{2}
    \label{Eqn:Posing_DB_BoLM_alternative} 
    &\frac{\mu}{k(\mathbf{x})} \mathbf{v} + \mathrm{grad}[p] - \mathrm{div}\big[2\mu \mathbf{D}[\mathbf{v}]\big]= \rho \mathbf{b}(\mathbf{x}) 
    &&\quad \mathrm{in} \; \Omega \\
    \label{Eqn:Posing_DB_Continuity_alternative} 
    &\mathrm{div}[\mathbf{v}] = 0 
    &&\quad \mathrm{in} \; \Omega \\
    \label{Eqn:Posing_DB_tBC_alternative} 
    &\big(-p \mathbf{I} + 2 \mu \mathbf{D}[\mathbf{v}]\big)\widehat{\mathbf{n}}(\mathbf{x}) = \mathbf{t}^{\mathrm{p}}(\mathbf{x}) 
    &&\quad \mathrm{on} \; \Gamma^{p} \\
    \label{Eqn:Posing_DB_vBC_alternative} 
    &\mathbf{v} = \mathbf{v}^{\mathrm{p}}(\mathbf{x}) 
    &&\quad \mathrm{on} \; \Gamma^{v} 
\end{alignat}
\end{subequations}
where $\mathbf{t}^{\mathrm{p}}(\mathbf{x})$ is the prescribed traction, which need not be in the form of a pressure. Here $\Gamma^{p}$ is the part of the boundary on which traction boundary condition is prescribed. The solution to above forward problem is again denoted by $\mathbf{v}_{\mathrm{DB}}(\mathbf{x};k(\mathbf{x}))$ and $p_{\mathrm{DB}}(\mathbf{x};k(\mathbf{x}))$.

Using the adjoint method, the sensitivity of the total dissipation rate to the permeability field can be written as follows: 
%---------------------------------------------------;
%  Equation: Sensitivities based on adjoint method  ;
%---------------------------------------------------;
{\small 
\begin{align}
D \Phi_{\mathrm{DB}}[k(\mathbf{x})] 
&= \int_{\Omega} \left(\frac{\mu}{k(\mathbf{x})} \mathbf{v}_{\mathrm{DB}}\big(\mathbf{x};k(\mathbf{x})\big) \cdot \mathbf{v}_{\mathrm{DB}}\big(\mathbf{x};k(\mathbf{x})\big) 
+ 2 \mu \,  \mathbf{D}[\mathbf{v}_{\mathrm{DB}}(\mathbf{x};k(\mathbf{x}))] \cdot 
\mathbf{D}[\mathbf{v}_{\mathrm{DB}}(\mathbf{x};k(\mathbf{x}))] 
\right)^{\prime} \, \mathrm{d} \Omega \notag \\
&\quad \quad -2 \int_{\Omega} \boldsymbol{\Lambda}(\mathbf{x}) \cdot \Big\{\underbrace{\frac{\mu}{k(\mathbf{x})} \mathbf{v}_{\mathrm{DB}}(\mathbf{x};k(\mathbf{x})) 
+ \mathrm{grad}[p_{\mathrm{DB}}(\mathbf{x};k(\mathbf{x}))] - \mathrm{div}\big[2 \mu \,  \mathbf{D}[\mathbf{v}_{\mathrm{DB}}(\mathbf{x};k(\mathbf{x}))]\big]- \rho \mathbf{b}(\mathbf{x})}_{Eq.~\eqref{Eqn:Posing_DB_BoLM_alternative}}\Big\}^{\prime} \mathrm{d} \Omega \notag \\
&\quad \quad +2 \int_{\Omega} 
\lambda(\mathbf{x}) \,  
\Big\{\underbrace{\mathrm{div}[\mathbf{v}_{\mathrm{DB}}(\mathbf{x};k(\mathbf{x}))]}_{Eq.~\eqref{Eqn:Posing_DB_Continuity_alternative}}\Big\}^{\prime} \, \mathrm{d} \Omega \notag \\
&\quad \quad -2 \int_{\Gamma^{p}} \boldsymbol{\Lambda}(\mathbf{x}) \cdot  \Big\{\underbrace{\Big(-p_{\mathrm{DB}}(\mathbf{x};k(\mathbf{x}))\mathbf{I} + 2\mu \mathbf{D}[\mathbf{v}_{\mathrm{DB}}(\mathbf{x};k(\mathbf{x}))]\Big)\widehat{\mathbf{n}}(\mathbf{x}) - \mathbf{t}^{\mathrm{p}}(\mathbf{x})}_{Eq.~\eqref{Eqn:Posing_DB_tBC_alternative}}\Big\}^{\prime} 
\mathrm{d} \Gamma \notag \\ 
&\quad \quad +2 \int_{\Gamma^{v}} \Big(-\lambda(\mathbf{x}) \mathbf{I} + 2 \mu \mathbf{D}[\boldsymbol{\Lambda}(\mathbf{x})]
- 2 \mu \mathbf{D}[\mathbf{v}_{\mathrm{DB}}(\mathbf{x};k(\mathbf{x}))]\Big)  \widehat{\mathbf{n}}(\mathbf{x}) \cdot \Big\{\underbrace{\mathbf{v}_{\mathrm{DB}}(\mathbf{x};k(\mathbf{x})) - \mathbf{v}^{\mathrm{p}}(\mathbf{x})}_{Eq.~\eqref{Eqn:Posing_DB_vBC_alternative}}\Big\}^{\prime} \mathrm{d} \Gamma
\end{align}}
By following a similar procedure as before, we get the following  adjoint problem corresponding to the alternative form of Darcy-Brinkman equations: 
%------------------------------------;
%  Equation: Adjoint problem for DB  ;
%------------------------------------;
{\small 
\begin{subequations}
\begin{alignat}{2}
    \label{Eqn:Posing_Adjoint_problem_DB_alternative_BoLM}
    &\frac{\mu}{k(\mathbf{x})} \boldsymbol{\Lambda}(\mathbf{x}) 
    + \mathrm{grad}[\lambda(\mathbf{x})] 
    - \mathrm{div}\big[2\mu   \mathbf{D}[\boldsymbol{\Lambda}(\mathbf{x})]\big] = \frac{\mu}{k(\mathbf{x})} \mathbf{v}_{\mathrm{DB}}\big(\mathbf{x};k(\mathbf{x})\big) 
    - \mathrm{div}\big[2 \mu  \mathbf{D}[\mathbf{v}_{\mathrm{DB}}(\mathbf{x};k(\mathbf{x}))]\big]
    && \; \mathrm{in} \; \Omega \\
    &\mathrm{div}[\boldsymbol{\Lambda}(\mathbf{x})] = 0 
    && \; \mathrm{in} \; \Omega \\
    &\Big(-\lambda(\mathbf{x}) \mathbf{I} + 2 \mu \, 
    \mathbf{D}[\boldsymbol{\Lambda}(\mathbf{x})]
    \Big) \widehat{\mathbf{n}}(\mathbf{x}) = 2 \mu \,  \mathbf{D}[\mathbf{v}_{\mathrm{DB}}(\mathbf{x};k(\mathbf{x}))] \widehat{\mathbf{n}}(\mathbf{x}) 
    && \; \mathrm{on} \; \Gamma^{p} \\
    \label{Eqn:Posing_Adjoint_problem_DB_alternative_vBC}
    &\boldsymbol{\Lambda}(\mathbf{x}) = \mathbf{0} 
    && \; \mathrm{on} \; \Gamma^{v} 
\end{alignat}
\end{subequations}
Even under this alternative form, 
the expression for the sensitivity of the total dissipation rate with respect to the permeability field is same as before:
%----------------------------------------;
%  Equation: Final sensitivities for DB  ;
%----------------------------------------;
\begin{align}
D \Phi_{\mathrm{DB}}[k(\mathbf{x})] &= \int_{\Omega} \frac{\mu}{k^2(\mathbf{x})} 
\mathbf{v}_{\mathrm{DB}}\big(\mathbf{x};k(\mathbf{x})\big) \cdot 
\Big(2 \boldsymbol{\Lambda}(\mathbf{x}) - \mathbf{v}_{\mathrm{DB}}\big(\mathbf{x};k(\mathbf{x})\big)\Big) \, \mathrm{d} \Omega 
\end{align}
where $\boldsymbol{\Lambda}(\mathbf{x})$ satisfies the adjoint problem \eqref{Eqn:Posing_Adjoint_problem_DB_alternative_BoLM}--\eqref{Eqn:Posing_Adjoint_problem_DB_alternative_vBC}. Due to a change in one of the boundary conditions, the definitions for the classes of boundary value problems should be modified accordingly for the alternative form of Darcy-Brinkman equations. 
%-------------------------------------;
%  Enumerate: Classification of BVPs  ;
%-------------------------------------;
\begin{enumerate}
\item \textbf{Class A:} Prescription of a traction on the entire boundary. Stated mathematically, 
\begin{align}
\Gamma^{v} = \emptyset
\end{align}
\item \textbf{Class B:} Prescription of a traction on a part of the boundary while homogeneous velocity boundary conditions on the rest. That is, 
\begin{align}
\Gamma^{p} \neq \emptyset, \; 
\Gamma^{v} \neq \emptyset, 
\quad \mathrm{and} \quad  \mathbf{v}^{\mathrm{p}}(\mathbf{x}) = \mathbf{0} \; \mathrm{on} \;  \Gamma^{v}
\end{align}
\item \textbf{Class C:} Prescription of a compatible velocity boundary condition on the entire boundary and requiring the prescribed body force to be a conservative vector field. Mathematically, 
\begin{align}
\Gamma^{p} = \emptyset, \; \int_{\Gamma^{v} = \partial \Omega} \mathbf{v}^{\mathrm{p}}(\mathbf{x}) \cdot \widehat{\mathbf{n}}(\mathbf{x}) \, \mathrm{d} \Gamma = 0 \quad \mathrm{and} \quad \rho\mathbf{b}(\mathbf{x}) = - \mathrm{grad}[\psi] 
\end{align}
where $\psi(\mathbf{x})$ is a scalar field. 
\item \textbf{Class D:} Prescription of velocity boundary conditions on a part of the boundary with zero traction on the rest, and requiring the body force is zero. Written mathematically, 
\begin{align}
\Gamma^{p} \neq \emptyset, \; \Gamma^{v} \neq 0, \;  \mathbf{t}^{\mathrm{p}}(\mathbf{x}) = \mathbf{0} \; \mathrm{on} \; \Gamma^{p}, \; \mathrm{and} \;  \mathbf{b}(\mathbf{x}) = \mathbf{0} 
\end{align}
\end{enumerate}

The results---solutions to the adjoint problem and expressions for the sensitivities---and conclusions remain the same as the original form of Darcy-Brinkman equations. Since the derivations are a straightforward extension of the original case, we omitted the details for brevity.

%================;
%  Bibliography  ;
%================;
\bibliographystyle{plainnat}
\bibliography{Master_References}
\end{document}